\documentclass[10pt,a4paper]{article}
\usepackage{bbm}

\usepackage[leqno]{amsmath}
\usepackage{amsfonts}
\usepackage{graphicx}
\usepackage{amsmath}
\usepackage{amssymb}
\usepackage{latexsym}
\usepackage{amsmath, amsfonts,amssymb, amsthm, euscript,makeidx,color,mathrsfs}
\usepackage{enumerate}

\usepackage[colorlinks,linkcolor=blue,anchorcolor=green,citecolor=red]{hyperref}%when print, use the next package
\def\sqr#1#2{{\vcenter{\vbox{\hrule height.#2pt
              \hbox{\vrule width.#2pt height#1pt \kern#1pt \vrule width.#2pt}
              \hrule height.#2pt}}}}
\def\signed #1{{\unskip\nobreak\hfil\penalty50
              \hskip2em\hbox{}\nobreak\hfil#1
              \parfillskip=0pt \finalhyphendemerits=0 \par}}
\def\endpf{\signed {$\sqr69$}}
\def\5n{\negthinspace \negthinspace \negthinspace \negthinspace \negthinspace }
\def\4n{\negthinspace \negthinspace \negthinspace \negthinspace }
\def\3n{\negthinspace \negthinspace \negthinspace }
\def\2n{\negthinspace \negthinspace }
\def\1n{\negthinspace }

\def\dbR{\mathbb{R}}
\def\dbS{\mathbb{S}}

\def\sK{\mathscr{K}}
\def\sL{\mathscr{L}}
\def\sM{\mathscr{M}}

\def\sR{\mathscr{R}}

\def\sU{\mathscr{U}}

\def\sX{\mathscr{X}}

\def\={\buildrel \triangle \over =}

\def\ds{\displaystyle}

\def\ns{\noalign{\ss}}
\def\b{\beta}

\def\d{\delta}
\def\e{\varepsilon}
\def\z{\zeta}
\def\k{\kappa}
\def\l{\lambda}

\def\si{\sigma}
\def\t{\tau}
\def\f{\varphi}
\def\th{\theta}
\def\o{\omega}

\def\i{\infty}
\def\G{\Gamma}
\def\D{\Delta}
\def\Th{\Theta}

\def\F{\Phi}

\def\cB{{\cal B}}

\def\cK{{\cal K}}

\def\BG{{\bf G}}

\def\BQ{{\bf Q}}
\def\BR{{\bf R}}
\def\BS{{\bf S}}

\def\Bg{{\bf g}}

\def\Bq{{\bf q}}

\def\BBG{\boldsymbol\Gamma}

\def\BBL{\boldsymbol\Lambda}

\def\BBl{\boldsymbol\lambda}

\def\BBrho{\boldsymbol\rho}

\def\no{\noindent}

\def\ss{\smallskip}
\def\ms{\medskip}
\def\bs{\bigskip}
\def\q{\quad}
\def\qq{\qquad}

\def\lan{\langle}
\def\ran{\rangle}

\def\rf{\eqref}
\def\h{\widehat}
\def\wt{\widetilde}

\def\cd{\cdot}
\def\cds{\cdots}

\def\ae{\hbox{\rm a.e.}}

\def\les{\leqslant}
\def\ges{\geqslant}

\def\({\Big (}
\def\){\Big )}
\def\[{\Big[}
\def\]{\Big]}

\def\bde{\begin{definition}\label}
\def\ede{\end{definition}}
\def\be{\begin{equation}}
\def\bel{\begin{equation}\label}
\def\ee{\end{equation}}
\def\bt{\begin{theorem}\label}
\def\et{\end{theorem}}
\def\bc{\begin{corollary}\label}
\def\ec{\end{corollary}}
\def\bl{\begin{lemma}\label}
\def\el{\end{lemma}}
\def\bp{\begin{proposition}\label}
\def\ep{\end{proposition}}
\def\bas{\begin{assumption}\label}
\def\eas{\end{assumption}}
\def\br{\begin{remark}\label}
\def\er{\end{remark}}
\def\bex{\begin{example}\label}
\def\ex{\end{example}}
\def\ba{\begin{array}}
\def\ea{\end{array}}
\def\ed{\end{document}}
\def\square#1{\vbox{\hrule\hbox{\vrule height#1%
     \kern#1\vrule}\hrule}}
\def\rectangle#1#2{\vbox{\hrule\hbox{\vrule height#1%
     \kern#2\vrule}\hrule}}

\font\tenbb=msbm10 \font\sevenbb=msbm7 \font\fivebb=msbm5

\newfam\bbfam
\scriptscriptfont\bbfam=\fivebb \textfont\bbfam=\tenbb
\scriptfont\bbfam=\sevenbb

\newtheorem{theorem}{\hskip 1.3em Theorem}[section]
\newtheorem{definition}[theorem]{\hskip 1.3em Definition}
\newtheorem{proposition}[theorem]{\hskip 1.3em Proposition}
\newtheorem{corollary}[theorem]{\hskip 1.3em Corollary}
\newtheorem{lemma}[theorem]{\hskip 1.3em Lemma}
\newtheorem{remark}[theorem]{\hskip 1.3em Remark}
\newtheorem{example}[theorem]{\hskip 1.3em Example}

\newtheorem{assumption}[theorem]{\hskip 1.3em Assumption}
\makeatletter
   
   \@addtoreset{equation}{section}
\makeatother

\begin{document}

\title{\bf Causal State Feedback Representation\\ for Linear Quadratic Optimal Control Problems of Singular Volterra Integral Equations
\footnote{This work was partially supported by the National Natural
Science Foundation of China under grant 12071067,
 National Key R\&D Program of China under grant 2020YFA0714102, and NSF grant DMS--1812921.}}
\date{}
\author{Shuo Han\footnote{  School of
Mathematics and Statistics, Northeast Normal University, Changchun
130024, China. E-mail: {\tt
hans861@nenu.} {\tt edu.cn.}},~~~ Ping Lin\footnote{  School of
Mathematics and Statistics, Northeast Normal University, Changchun
130024, China. E-mail: {\tt linp258@nenu.} {\tt edu.cn}.}~~~~and~~~~Jiongmin Yong \footnote{Department of Mathematics, University of Central Florida, Orlando, FL 32816, USA.  E-mail: {\tt jiongmin.yong@ucf.edu.}}}

\ms

\maketitle

{\bf Abstract.} This paper is concerned with a linear quadratic optimal control
for a class of singular Volterra integral equations. Under proper convexity
conditions, optimal control uniquely exists, and it could be characterized via Fr\'echet derivative of the quadratic functional in a Hilbert space or via maximum principle type necessary conditions. However, these (equivalent) characterizations have a shortcoming that the current value of the optimal control depends on the future values of the optimal state. Practically, this is not feasible. The main purpose of this paper is to obtain a causal state feedback representation of the optimal control.

\bs

\ms


\section{Introduction.}

Consider the following controlled singular linear Volterra integral equation:
\bel{state}X(t)=\f(t)+\int_0^t{A(t,s)X(s)+B(t,s)u(s)\over (t-s)^{1-\b}}ds,\qq \ae~t\in[0,T].\ee
In the above, $T>0$ is a fixed finite time horizon, $\f(\cd)$ is a given map, called the {\it free term} of the state equation, $X(\cd)$ is called the {\it state trajectory} taking values in the Euclidean space $\dbR^n$, $u(\cd)$ is called the {\it control} taking values in the Euclidean space $\dbR^m$, $A(\cd\,,\cd)$ and $B(\cd\,,\cd)$ are called the coefficients, taking values in $\dbR^{n\times n}$ and $\dbR^{n\times m},$ respectively, and $\b>0$.

\ms

We denote $\sX=L^2(0,T;\dbR^n)$, $\sU=L^2(0,T;\dbR^m)$. Under some mild conditions, for any control $u(\cd)\in\sU$, the state equation \rf{state} admits a unique solution $X(\cd)\in \sX$. To measure the performance of the control, we introduce the following quadratic cost functional
\bel{cost}\ba{ll}
\ns\ds J(u(\cd))=\int_0^T\(\lan Q(t)X(t),X(t)\ran+2\lan S(t)X(t),u(t)\ran+\lan R(t)u(t),u(t)\ran\\
\ns\ds\qq\qq\qq\qq+2\lan q(t),X(t)\ran+2\lan\rho(t),u(t)\ran\)dt+\lan GX(T),X(T)\ran+2\lan g,X(T)\ran,\ea\ee
where $Q(\cd)\in L^\infty(0,T;\dbS^n)$, $S(\cd)\in L^\infty(0,T;\dbR^{m\times n})$, $R(\cd)\in L^\infty(0,T;\dbS^m)$, $q(\cd)\in \sX$, $\rho(\cd)\in \sU$, $G\in\dbS^n$, $g\in\dbR^n$, with $\dbS^k$ being the set of all $(k\times k)$ symmetric (real) matrices. Our optimal control problem can be stated as follows.

\ms

\no\bf Problem (P). \rm Find a control $\bar u(\cd)\in\sU$ such that
\bel{OPC}J(\bar u(\cd))=\inf_{u(\cd)\in\sU}J(u(\cd)).\ee
Any $\bar u(\cd)$ satisfying \eqref{OPC} is called an {\it open-loop optimal control} of Problem (P), the corresponding state $\bar X(\cd)$ is called an {\it open-loop optimal state} and $(\bar X(\cd),\bar u(\cd))$ is called an {\it open-loop optimal pair}.

\ms

Memory exists in many application problems, heat transfer, population growth, disease spread, to mention a few. Volterra integral equations can be used to describe some dynamics involving memories. Study of optimal control problems for Volterra integral equations can be traced back to the works of Vinokurov in the later 1960s \cite{Vinokurov 1969}, followed by the works of Angell \cite{Angell 1976}, Kamien-Muller \cite{Kamien-Muller 1976}, Medhin \cite{Medhin 1986}, Carlson \cite{Carlson 1987}, Burnap--Kazemi \cite{Burnap-Kazemi 1999}, and some recent works by de la Vega \cite{de la Vega 2006}, Belbas
\cite{Belbas 2007,Belbas 2008}, and Bonnans--de la Vega--Dupuis \cite{Bonnans-de la Vega-Dupuis 2013}. All of the above-mentioned works are concerned with non-singular Volterra integral equations which exclude the case of \rf{state} with $\b\in(0,1)$. On the other hand, in the past several decades, fractional (order) differential equations have attracted quite a few researchers' attention due to some very interesting applications in physics, chemistry, engineering, population dynamics, finance and other sciences; See Oldham--Spanier \cite{Oldham-Spanier 1974} for some early examples of diffusion processes, Torvik--Bagley \cite{Torvik-Bagley 1984}, Caputo \cite{Caputo 1967}, and Caputo--Mainardi \cite{Caputo--Mainardi 1971} for modeling of the mechanical properties of materials, Benson \cite{Benson 1998} for the advection and the dispersion of solutes in natural porous or fractured media, Chern \cite{Chern 1993}, Diethelm--Freed \cite{Diethelm-Freed 1999} for the modeling  behavior of viscoelastic and viscoplastic materials under external influences,
Scalas--Gorenflo--Mainardi \cite{Scalas-Gorenflo-Mainardi 2004} for the mathematical models in finance, Das--Gupta \cite{Das-Gupta 2011}, Demirci--Unal--\"Ozalp \cite{Demirci-Unal-Ozalp 2011}, Arafa--Rida--Khalil \cite{Arafa-Rida-Khalil 2012}, Diethelm \cite{Diethelm 2013} for some population and epidemic models, Metzler et al. \cite{Metzler-Schick-Kilian-Nonnenmacher 1995} for the relaxation in filled polymer networks, and Okyere et al. \cite{Okyere-Oduro-Amponsah-Dontwi-Fempong 2016} for a SIR model with constant population. An extensive survey on fractional differential equations can be found in the book by Kilbas--Srivastava--Trujillo \cite{Kilbas-Srivastava-Trujillo 2006}. In the recent years, optimal control problems have been studied for fractional differential equations by a number of authors. We mention the works of Agrawal \cite{Agrawal 2004, Agrawal 2008}, Agrawal--Defterli--Baleanu \cite{Agrawal-Defterli-Baleanu 2010},  Bourdin \cite{Bourdin 2012}, Frederico--Torres \cite{Frederico-Torres 2008}, Hasan--Tangpong--Agrawal \cite{Hasan-Tangpong-Agrawal} and Kamocki \cite{Kamocki 2014,Kamocki 2014b}, Gomoyunov \cite{Gomoyunov 2020}, Koenig \cite{Koenig 2020}.

\ms

It turns out that fractional differential equations (of the order no more than 1), in the sense of Riemann--Liouville or in the sense of Caputo, are equivalent to Volterra integral equations with the integrand being singular along $s=t$, and the free term $\f(\cd)$ being possibly discontinuous (blowing up) at $t=0$ (See \cite{Lin-Yong 2020} for some details). More precisely, in the linear case, the corresponding controlled state equation of form \rf{state} could have the free term look like the following:
\bel{1.3*}\f(t)={c\over t^{1-\b}}\ (\mbox{or}\ c),\ee
for some constant $c\in\dbR$. In \cite{Lin-Yong 2020}, a class of controlled nonlinear singular Volterra integral equations was considered. Well-posedness of the state equation and some regularity of the state trajectory were established, and a Pontryagin type maximum principle for optimal controls was proved.

\ms

On the other hand, Pritchard--You \cite{Prichard-You 1996} considered the quadratic optimal control problems for the following controlled linear Volterra integral equations in a Hilbert space $H$:
\bel{K}y(t)=f(t)+\int_0^tF(t,\t)u(\t)d\t,\qq t\in[0,T].\ee
It was assumed in \cite{Prichard-You 1996} that $f(\cd)\in C([0,T];H)$ and $F:\bar\D\to\sL(U;H)$ is strongly continuous in the sense that for each $u\in U$, $F(\cd\,,\cd)u\in C(\bar\D;H)$. Here, $U$ is another Hilbert space and
$\bar\D$ is the closure of the following set
\bel{D}\D=\{(t,s)\in[0,T]^2\bigm|0\les s<t\les T]\}.\ee
Thus, in particular, the following holds
\bel{K<infty}\|F(t,t)\|_{\sL(U;H)}<\infty,\qq t\in[0,T].\ee
This excludes our state equation \rf{state} which has a singular kernel. We will see later that when the variation of constants formula is applied, our state process will have a similar representation as \rf{K}, but with both the free term $f(\cd)$ and the operator $F(\cd\,,\cd)$ being not necessarily continuous.

\ms

Practically, if an optimal control exists, one expects that the optimal could have a state feedback representation which is non-anticipating. In the case of state equation being an ordinary differential equation (or a partial differential equation, a stochastic differential equation), such kind of representation can be obtained, under some proper conditions, via a solution to a Riccati differential equation.
Pandolfi \cite{Pandolfi 2018} derived an optimal feedback control for a Volterra integro-differential equation by using the corresponding Riccati equation. That could be done because the state equation in \cite{Pandolfi 2018} was of a special form which has the semigroup property and thus one could use semigroup representation to derive a theory of Riccati equation in a standard way. But the general Volterra equation does not have a semigroup evolutionary property. For the controlled linear Volterra integral equation of form \rf{K}, a so-called {\it projection causality} approach was introduced in \cite{Prichard-You 1996}. The optimal control could be represented as a so-called linear {\it causal feedback} (see later for a precise definition) of the state trajectory with the feedback operator being determined by a solution to a Fredholm integral equation.

\ms

In this paper, we will carry out some careful analysis on the state equation, and pay special attention to certain continuity of the state trajectory since in the quadratic cost functional, the terminal value $X(T)$ of the state trajectory is involved. Also, we make it clear that the condition $\b>{1\over2}$ should be assumed in order the LQ problem is well-formulated. By a standard method for minimization of a quadratic functional in Hilbert space, we obtain a characterization of the (open-loop) optimal control in some abstract form. On the other hand, by variational method, in the spirit of maximum principle, we may obtain another characterization of the open-loop optimal control. We will show that these two characterizations are equivalent. However, from those characterizations, the open-loop optimal control is not non-anticipating in the sense that in determining the value $\bar u(t)$ of the open-loop optimal control $\bar u(\cd)$ at time $t$, some future information $\{\bar X(s)\bigm|s\in[t,T]\}$ of the optimal state trajectory has to be used. This is not practically realizable. In the classical LQ problem of differential equations, one could get a closed-loop representation of the open-loop optimal control via the solution to a Riccati equation. However, for general integral equations, such an approach is not working. In fact, the problem we considered in this paper is a nonlocal problem. Thus, we could not obtain a closed-loop optimal control whose current value only depends on the current state value, by using the standard method of ODE (or PDE, Volterra integro-differential equation in \cite{Pandolfi 2018}) in terms of the Riccati equation.
Inspired by \cite{Prichard-You 1996}, we will try to obtain a causal state feedback representation for the open-loop optimal control in the following sense: The current value $\bar u(t)$ of the open-loop optimal control $\bar u(\cd)$ is written in terms of the current optimal state value $\bar X(t)$, as well as a {\it causal trajectory} $\bar X_t(\cd)$ and an {\it auxiliary trajectory} $\bar X^a(t)$, via a family of Fredholm integral equations which essentially plays a role of Riccati equation in the classical LQ problems. It is worthy of pointing out that $X_t(\cd)$ and $X^a(\cd)$ can be running at the same time as the state equation, and they are non-anticipating. Note that $X_t(\cd)$ and $X^a(\cd)$ are not involved in calculating the cost functional, but they are used to represent the open-loop optimal control. Although the main idea comes from \cite{Prichard-You 1996}, our modified version of the method is more direct which reveals the essence of the problem more clearly.
In the proof of \cite{Prichard-You 1996}, they introduce an abstract operator to establish the interrelations between the state trajectory and the causal trajectory. In this paper, we do not need to introduce the similar abstract operator, but give a more direct proof. Furthermore, the trajectory $X_t(T)$ in \cite{Prichard-You 1996} may raise the doubt about the causality. In this paper, we can avoid this doubt by introducing the auxiliary trajectory $X^a(\cd)$.
\ms

The rest of the paper is organized as follows. In section 2, we carry out some
analysis for the state equation. Section 3 is devoted to the open-loop optimal control and its characterizations. Causal projection as well as abstract form of casual state feedback representation of the open-loop optimal control is presented in Section 4. We introduce a family of Fredholm integral equations in Section 5, which makes the representation obtained in Section 4 more practically accessible. In section 6, we briefly present a possible numerical scheme which is applicable to solve the Fredholm integral equation obtained in Section 5.

\section{Preliminary Results.}

In this section, we will present some preliminary results which will be useful later. Let us recall $\D$ defined by \rf{D}. Note that the ``diagonal line'' $\{(t,t)\,|\,t\in[0,T]\}$ is not contained in $\D$. Thus if $(t,s)\mapsto f(t,s)$ is continuous on $\D$, $f(\cd\,,\cd)$ is allowed to be unbounded as $|t-s|\to0$. Throughout this paper, we denote $t_1\vee t_2=\max\{t_1,t_2\}$ and $t_1\land t_2=\min\{t_1,t_2\}$, for any $t_1,t_2\in\dbR$. The characteristic function of any set $E$ is denoted by ${\bf1}_E(\cd)$. For any set $E\subseteq\dbR$ and a function $\f:E\to\dbR$, we extend it to be zero in $\dbR\setminus E$. We call a strictly increasing continuous function $\o(\cd):[0,\infty)\to[0,\infty)$ a {\it modulus of continuity} if $\o(0)=0$. Also, $K$ will be a generic constant which could be different from line to line.

\ms

Let us recall the Young's inequality for convolution (Theorem 3.9.4 in \cite{Bogachev 2007}).

\bl{Young} \sl Let $p,q,r\in[1,+\infty]$ satisfy ${1\over p}+1={1\over q}+{1\over r}.$
Then for any $f(\cd)\in L^q(\dbR^n)$, $g(\cd)\in L^r(\dbR^n)$,
\bel{}\|f(\cd)*g(\cd)\|_{L^p(\dbR^n)}\les\|f(\cd)\|_{L^q(\dbR^n)}\|g(\cd)\|_{L^r(\dbR^n)}.\ee

\el

From the above lemma, we have the following corollary which is a refinement of
that found in \cite{Lin-Yong 2020}.

\bc{corollary 2.2} \sl Let $\th:\D\to\dbR^n$ and $\th_0:[0,T]\to\dbR$ be measurable
such that
\bel{|th_0|}|\th(t,\t)|\les\th_0(\t),\qq\ae~(t,\t)\in\D.\ee
For any $s\in[0,T)$, define
$$\eta(t,s)=\int_s^t{\th(t,\t)\over(t-\t)^{1-\b}}d\t,\qq t\in(s,T].$$

{\rm(i)} Let $\b\in(0,1)$, $1\les r<{1\over1-\b}$, ${1\over p}+1={1\over q}+{1\over r}$, $p,q\in[1,\infty]$, and $\th_0(\cd)\in L^q(s,T)$.
Then
\bel{|eta|}\|\eta(\cd,s)\|_{L^p(s,T;\dbR^n)}\les\({(T-s)^{1-r(1-\b)}\over1-r(1-\b)}\)^{
1\over r}\|\th_0(\cd)\|_{L^q(s,T)}.\ee
In particular, with $q=2$,
\bel{|eta|}\|\eta(\cd,s)\|_{L^p(s,T;\dbR^n)}\les{(T-s)^{\b-({1\over2}-{1\over p})}\over\({\b-({1\over2}-{1\over p})\over{1\over2}+{1\over p}}\)^{{1\over2}+{1\over p}}}\|\th_0(\cd)\|_{L^2(s,T)}.\ee
Further, with both $p=q=2$,
\bel{|eta|2}\|\eta(\cd,s)\|_{L^2(s,T;\dbR^n)}\les{(T-s)^\b\over\b}\|\th_0(\cd)
\|_{L^2(s,T)},\ee
and with $p=\i$, $q=2$, $\b>{1\over2}$,
\bel{|eta|infty}\|\eta(\cd,s)\|_{L^\infty(s,T;\dbR^n)}\les{(T-s)^{\b-{1\over2}}
\over\sqrt{2\b-1}}\|\th_0(\cd)\|_{L^2(s,T)}.\ee

{\rm(ii)} Let there exist a $\d_0\in(0,\frac{T-s}{2}),$ and a modulus of continuity $\o$ such that $\forall t\in[T-2\delta_0,T],$
\bel{}\ba{ll}
\ns\ds|\th(t,\t)-\th(T,\t)|\les\o(T-t)\th_1(\t),\qq\ae~ \t\in[s,t),\\
\ns\ds|\th(t,\t)|\les\th_0(\t),\qq\ae~ \t\in[s,t),\ea\ee
for some $\th_0(\cd)$, $\th_1(\cd)\in L^q(s,T)$, $q\in(1,\infty]$, $1>\b>{1\over q}$.
Then $\eta(\cd,s)$ is continuous at $T$.

\ec

\it Proof. \rm (i) By Lemma \ref{Young} with $f(\t)=\th_0(\t){\bf1}_{[s,T]}(\tau)$ and $g(\t)
={1\over\t^{1-\b}}{\bf1}_{(0,T-s]}(\tau)$, we can obtain our conclusion. The rest is clear.

\ms

(ii) Pick any $\d\in(0,\d_0)$. For any $t\in(T-\d,T)$, we look at the following,
assuming first that $q\in(1,\infty)$ and setting $\k=(1-\b){q\over q-1}<1$ (since
$\b>{1\over q}$):
$$\ba{ll}
\ns\ds|\eta(T,s)-\eta(t,s)|=\Big|\int_s^T{\th(T,\t)\over(T-\t)^{1-\b}}d\t
-\int_s^t{\th(t,\t)\over(t-\t)^{1-\b}}d\t\Big|\\
\ns\ds\les\int_s^{t-\d}|\th(t,\t)|\({1\over(t-\t)^{1-\b}}-{1\over(T-\t)^{1-\b}}\)d\t
+\int_s^{t-\d}{|\th(t,\t)-\th(T,\t)|\over(T-\t)^{1-\b}}d\t\\
\ns\ds\qq+\int_{t-\d}^t{\th_0(\t)\over(t-\t)^{1-\b}}d\t+\int_{t-\d}^T
{\th_0(\t)\over(T-\t)^{1-\b}}d\t\\
\ns\ds\les(T-t)^{1-\b}\int_s^{t-\d}{\th_0(\t)\over(t-\t)^{1-\b}
(T-\t)^{1-\b}}d\t+\o(T-t)\int_s^{t-\d}{\th_1(\t)\over(T-\t)^{1-\b}}d\t\\
\ns\ds\qq+\|\th_0(\cd)\|_{L^q(s,T)}\(\int_{t-\d}^t{d\t\over(t-\t)^\k}
\)^{q-1\over q}+\|\th_0(\cd)\|_{L^q(s,T)}\(\int_{t-\d}^T{d\t\over(T-\t)^\k}
\)^{q-1\over q}\\
\ns\ds\les{(T-t)^{1-\b}\over{\d^{2(1-\b)}}}\|\th_0(\cd)\|_{L^1(s,T)}
+\o(T-t)\({(T-s)^{1-\k}\over1-\k}\)^{q-1\over q}\|\th_1(\cd)\|_{L^q(s,T)}\ea$$

$$\ba{ll}
\ns\ds+\|\th_0(\cd)\|_{L^q(s,T)}\[\({\d^{1-\k}\over1-\k}\)^{q-1\over q}+\({(T-t+\d)^{1-\k}\over1-\k}\)^{q-1\over q}\].\ea$$
Hence, for any $\e>0$, we first take $\d>0$ sufficiently small so that
$$\|\th_0(\cd)\|_{L^q(s,T)}\[\({\d^{1-\k}\over1-\k}\)^{q-1\over q}+\({(2\d)^{1-\k}\over1-\k}\)^{q-1\over q}\]<{\e\over2}.$$
Since the modulus of continuity $\omega(\cd)$ is continuous and $\omega(0)=0$, we  can take $\bar\d\in(0,\d)$ even smaller so that
$${\bar\d^{1-\b}\over{\d^{2(1-\b)}}}\|\th_0(\cd)\|_{L^1(s,T)}+\o(\bar\d)
\({(T-s)^{1-\k}\over1-\k}\)^{q-1\over q}\|\th_1(\cd)\|_{L^q(s,T)}<{\e\over2}.$$
Combining the above, we see that $\eta(\cd,s)$ is continuous at $T$.

\ms

In the case $q=\infty$, we have
$$\ba{ll}
\ns\ds|\eta(T,s)-\eta(t,s)|=\Big|\int_s^T{\th(T,\t)\over(T-\t)^{1-\b}}d\t
-\int_s^t{\th(t,\t)\over(t-\t)^{1-\b}}d\t\Big|\\
\ns\ds\les{(T-t)^{1-\b}\over{\d^{2(1-\b)}}}\|\th_0(\cd)\|_{L^1(s,T)}
+\o(T-t){(T-s)^\b\over\b}\|\th_1(\cd)\|_{L^\infty(s,T)}\\
\ns\ds\qq+\|\th_0(\cd)\|_{L^\infty(s,T)}\[{\d^\b\over\b}+{(T-t+\d)^\b\over\b}\].\ea$$
Then, similar to the above, we obtain the continuity of $\eta(\cd,s)$ at
$T$. \endpf

\ms

We now look at the following linear Volterra integral equation
\bel{Volterra1}X(t)=\xi(t)+\int_0^t{A(t,s)X(s)\over(t-s)^{1-\b}}ds,\qq \ae~t\in[0,T].\ee
Note that \rf{state} is a case of the above with
\bel{xi}\xi(t)=\f(t)+\int_0^t{B(t,s)u(s)\over(t-s)^{1-\b}}ds,\qq \ae~t\in[0,T].\ee
Before going further, we introduce the following assumption for the coefficients of \rf{state}.

\ms

{\bf(A1)} The coefficients $A(\cd\,,\cd)\in L^\infty(\D;\dbR^{n\times n})$ and $B(\cd\,,\cd)\in L^\infty(\D;\dbR^{n\times m})$. The free term $\f(\cd)\in\sX$.

\ms

For convenience, throughout the paper, we assume that
$$|A(t,s)|\les\|A\|_\infty,\q|B(t,s)|\les\|B\|_\infty,\qq\forall(t,s)\in\D.$$
Then, under (A1), for any $u(\cd)\in\sU$, by Corollary \ref{corollary 2.2},
(i), we have that if $\beta\in(0,1),$
\bel{Bu}\ba{ll}
\ns\ds\(\int_0^T\Big|\int_0^t{B(t,s)u(s)\over(t-s)^{1-\b}}ds\Big|^2dt\)^{1\over2}
\les\|B\|_\i\[\int_0^T
\(\int_0^t{|u(s)|\over(t-s)^{1-\b}}ds\)^2dt\]^{1\over2}\les\|B\|_\i\({T^\b\over\b}\)
\|u(\cd)\|_\sU.\ea\ee
Consequently, for the free term $\xi(\cd)$ defined by \rf{xi}, one has
\bel{|xi|}\|\xi(\cd)\|_\sX\les\|\f(\cd)\|_\sX+{T^\b\|B\|_\i\over\b}
\|u(\cd)\|_\sU.\ee
The following gives the well-posedness of \rf{Volterra1} as well as its variation of
constants formula.

\bt{well-posed} \sl Let {\rm(A1)} hold. Then for any $\xi(\cd)\in\sX$, \rf{Volterra1} admits a unique solution $X(\cd)\in\sX$. Moreover, there exists a measurable function
$\F(\cd\,,\cd):\D\to\dbR^{n\times n}$ satisfying that for any $s\in[0,T),$
\bel{Phi}\F(t,s)={A(t,s)\over(t-s)^{1-\b}}+\int_s^t{A(t,\t)\F(\t,s)\over (t-\t)^{1-\b}}d\t,\qq t\in(s,T],\ee
such that for some constant $K>0$,
\bel{|Phi|}|\F(t,s)|\les{K\over(t-s)^{1-\b}},\qq (t,s)\in\D,\ee
and the solution $X(\cd)$ to \rf{Volterra1} can be represented by the
following {\it variation of constants formula}:
\bel{variation}X(t)=\xi(t)+\int_0^t\F(t,s)\xi(s)ds,\qq\ae~ t\in[0,T],\ee
with the following estimate:
\bel{|X|}\|X(\cd)\|_\sX\les K\|\xi(\cd)\|_\sX.\ee
Moreover, the function $\F(\cd\,,\cd)$ also satisfies that for any $s\in[0,T),$
\bel{Phi*}\F(t,s)={A(t,s)\over(t-s)^{1-\b}}+\int_s^t{\F(t,\t)A(\t,s)\over (\t-s)^{1-\b}}d\t,\qq t\in(s,T].\ee

\et

\it Proof. \rm First of all, by a standard contraction mapping argument,
making use of Corollary \ref{corollary 2.2}, (i), we see that for any $\xi(\cd)\in\sX$, \rf{Volterra1} admits a unique solution $X(\cd)\in\sX$.

\ms

Next, by (A1), if $\F(\cd\,,\cd)$ is a
solution of \rf{Phi}, then
$$|\F(t,s)|\les{\|A\|_\i\over(t-s)^{1-\b}}+\int_s^t{\|A\|_\i|\F(\t,s)|
\over(t-\t)^{1-\b}}d\t,\qq(t,s)\in\D.$$
Thus, by Gronwall's inequality, we have
$$\ba{ll}
\ns\ds|\F(t,s)|\les{\|A\|_\i\over(t-s)^{1-\b}}+K\int_s^t{\|A\|^2_\i\over
(t-\t)^{1-\b}(\t-s)^{1-\b}}d\t\\
\ns\ds={\|A\|_\i\over(t-s)^{1-\b}}+{K\|A\|^2_\i\cB(\b,\b)\over(t-s)^{1-2\b}}
={\|A\|_\i+K\|A\|^2_\i\cB(\b,\b)(t-s)^\b\over(t-s)^{1-\b}}\les{K\over
(t-s)^{1-\b}},\qq(t,s)\in\D.\ea$$
This proves \rf{|Phi|}. In the above, $\cB(\cd\,,\cd)$ is the Beta function,
and recall that $K$ stands for a generic constant which could be different from line
to line.

\ms

Now, we inductively define the following sequence of measurable functions:
\bel{F_k}F_1(t,s)={A(t,s)\over(t-s)^{1-\b}},\qq F_{k+1}(t,s)=\int_s^tF_1(t,\t)
F_k(\t,s)d\t,\q k=1,2,3,\cds,\qq(t,s)\in\D.\ee
Then
$$\ba{ll}
\ns\ds|F_1(t,s)|\les{\|A\|_\i\over(t-s)^{1-\b}},\qq(t,s)\in\D,\\
\ns\ds|F_2(t,s)|\les\int_s^t|F_1(t,\t)F_1(\t,s)|d\t\les\|A\|_\i^2\int_s^t{d\t\over(t-\t)^{1-\b}
(\t-s)^{1-\b}}={\|A\|_\i^2\cB(\b,\b)\over(t-s)^{1-2\b}},\qq(t,s)\in\D.\ea$$
By induction, we can show that
\bel{|L_k|}|F_k(t,s)|\les\int_s^t|F_1(t,\t)F_{k-1}(\t,s)|d\t\les{\|A\|_\i^k\over(t-s)^{1-k\b}}
\prod_{j=1}^{k-1}\cB(\b,j\b),\q k\ges1,\qq(t,s)\in\D.\ee
According to \cite{Spiegel 1969} (p.102), Gamma function $\G(\cd)$ admits the following asymptotic expansion:
$$\G(z+1)=\sqrt{2\pi z}\,z^ze^{-z}\(1+R(z)\),\qq z>\3n>1;\qq R(z)={1\over12z}+{1\over288z^2}-{139\over51840z^3}+\cds.$$
Thus, for $j$ large enough, we have
$$\ba{ll}
\ns\ds\cB(\b,j\b)={\G(\b)\G(j\b)\over\G(j\b+\b)}={\G(\b)\sqrt{2\pi(j\b-1)}
(j\b-1)^{j\b-1}e^{-(j\b-1)}\(1+R(j\b-1)\)\over\sqrt{2\pi(j\b+\b-1)}(j\b+\b-1)^{
j\b+\b-1}e^{-(j\b+\b-1)}\(1+R(j\b+\b-1)\)}\\
\ns\ds=\G(\b)\sqrt{j\b-1\over j\b+\b-1}\({j\b-1\over j\b+\b-1}\)^{j\b-1}{e^\b\over
(j\b+\b-1)^\b}\(1+\wt R(j)\)\\
\ns\ds\les{\G(\b)e^\b\over\b^\b}\({j\over j+1-{1\over\b}}\)^\b\(1+\wt R(j)\){1\over j^\b}\les{\G(\b)e^\b\over\b^\b}{2^{\b+1}\over j^\b},\ea$$
for some $\wt R(j)\to0$ as $j\to\infty$. Consequently, there exists a $k_0$ such that for $k>k_0$,
$$\ba{ll}
\ns\ds|F_k(t,s)|\1n\les\1n{\|A\|_\i^k\over(t\1n-\1n s)^{1-k\b}}
\(\prod_{j=1}^{k_0-1}\cB(\b,j\b)\)\({\G(\b)
e^\b2^{\b+1}\over\b^\b}\)^{k-k_0}{[(k_0\1n-\1n1)!]^\b\over[(k\1n-\1n1)!]^\b}
\1n\les\1n{K_0K^k\over
(t\1n-\1n s)^{1-k\b}}{1\over[(k\1n-\1n1)!]^\b},\q(t,s)\1n\in\1n\D,\ea$$
for some constants $K_0,K>0$. Then, for any $\d>0$, series
$\ds\sum_{k=1}^\i F_k(t,s)$ is uniformly and absolutely convergent
for $(t,s)\in\D_\d$ with
$$\D_\d=\{(t,s)\in\D\bigm|t-s\ges\d\big\}.$$
Now we define
$$\F(t,s)=\sum_{k=1}^\i F_k(t,s),\qq(t,s)\in\D,$$
which is measurable (since each $F_k(\cd\,,\cd)$ is measurable) and
bounded on each $\D_\d$, $\d>0$. We can easily check that the above
defined $\F(\cd\,,\cd)$ is the unique solution of \rf{Phi}, and therefore,
estimate \rf{|Phi|} holds. Further, for any $\xi(\cd)\in L^2(0,T;\dbR^n)$,
similar to \rf{Bu}, we see that
$$X(t)=\xi(t)+\int_0^t\F(t,s)\xi(s)ds,\qq\ae~ t\in[0,T],$$
is well-defined as an element in $L^2(0,T;\dbR^n)$. In addition, for such
defined $X(\cd)$, one has
$$\ba{ll}
\ns\ds\int_0^t{A(t,s)X(s)\over(t-s)^{1-\b}}ds=\int_0^t{A(t,s)\over
(t-s)^{1-\b}}\[\xi(s)+\int_0^s\F(s,\t)\xi(\t)d\t\]ds\\
\ns\ds=\int_0^t{A(t,s)\xi(s)\over(t-s)^{1-\b}}ds+\int_0^t\int_\t^t
{A(t,s)\F(s,\t)\xi(\t)\over(t-s)^{1-\b}}dsd\t\\
\ns\ds=\int_0^t\({A(t,s)\over(t-s)^{1-\b}}+\int_s^t{A(t,\t)\F(\t,s)\over
(t-\t)^{1-\b}}d\t\)\xi(s)ds=\int_0^t\F(t,s)\xi(s)ds=X(t)-\xi(t),\q\ae~t\in[0,T].\ea$$
This proves \rf{variation}. Making use of \rf{|Phi|}, and similar to
\rf{Bu}, we obtain
$$\ba{ll}
\ns\ds\|X(\cd)\|_\sX\les\|\xi(\cd)\|_\sX+\(\int_0^T\Big|\int_0^t\F(t,s)\xi(s)ds
\Big|^2dt\)^{1\over2}\\
\ns\ds\qq\qq\les\|\xi(\cd)\|_\sX+\(\int_0^T\(\int_0^t{|\xi(s)|\over(t-s)^{1-\b}}ds
\)^2dt\)^{1\over2}\les\|\xi(\cd)\|_\sX+K\|\xi(\cd)\|_\sX.\ea$$
Finally, we prove \rf{Phi*}. To this end, we make the following observation:
$$\ba{ll}
\ns\ds F_3(t,s)=\int_s^tF_1(t,\t)F_2(\t,s)d\t=\int_s^tF_1(t,\t)\int_s^\t
F_1(\t,r)F_1(r,s)drd\t\ea$$
$$\ba{ll}
\ns\ds\qq\q=\int_s^t\int_r^tF_1(t,\t)F_1(\t,r)F_1(r,s)d\t dr
=\int_s^tF_2(t,r)F_1(r,s)dr,\qq(t,s)\in\D.\ea$$
Hence, by induction, we see that (comparing with \rf{F_k})
$$F_{k+1}(t,s)=\int_s^tF_k(t,r)F_1(r,s)dr, \qq(t,s)\in\D,\ ~k\ges1.$$
Then we see that \rf{Phi*} holds. This completes the proof. \endpf

\ms

According to the above theorem, for state equation \rf{state}, we have the
following representation of the state process $X(\cd)$ in terms of the
control $u(\cd)$ and the free term $\f(\cd)$:
\bel{X=u}\ba{ll}
\ns\ds X(t)=\f(t)+\int_0^t{B(t,s)u(s)\over(t-s)^{1-\b}}ds
+\int_0^t\F(t,s)\[\f(s)+\int_0^s{B(s,\t)u(\t)\over(s-\t)^{1-\b}}d\t\]ds\\
\ns\ds=\f(t)+\int_0^t\F(t,s)\f(s)ds+\int_0^t{B(t,s)u(s)\over(t-s)^{1-\b}}
ds+\int_0^t\int_\t^t{\F(t,s)B(s,\t)u(\t)\over(s-\t)^{1-\b}}dsd\t\\
\ns\ds=\f(t)+\int_0^t\F(t,s)\f(s)ds+\int_0^t\({B(t,s)\over(t-s)^{1-\b}}
+\int_s^t{\F(t,\t)B(\t,s)\over(\t-s)^{1-\b}}d\t\)u(s)ds\\
\ns\ds\equiv\psi(t)+\int_0^t\Psi(t,s)u(s)ds,\ea\ee
where
\bel{psi}\ba{ll}
\ns\ds\psi(t)=\f(t)+\int_0^t\F(t,s)\f(s)ds,\qq\ae~ t\in[0,T],\\
\ns\ds\Psi(t,s)={B(t,s)\over(t-s)^{1-\b}}
+\int_s^t{\F(t,\t)B(\t,s)\over(\t-s)^{1-\b}}d\t,\qq(t,s)\in\D.\ea\ee
Clearly, $\Psi:\D\to\dbR^{n\times m}$, and
\bel{|Psi|}\ba{ll}
\ns\ds|\Psi(t,s)|\les{\|B\|_\i\over(t-s)^{1-\b}}+\|B\|_\i\int_s^t{|\F(t,\t)|\over
(\t-s)^{1-\b}}d\t\\
\ns\ds\qq\q~\les{K\over(t-s)^{1-\b}}+K\int_s^t{d\t\over(t-\t)^{1-\b}(\t-s)^{1-\b}}
\les{K\over(t-s)^{1-\b}},\qq(t,s)\in\D.\ea\ee
Moreover, noting $\xi(\cd)$ defined by \rf{xi} and the estimate \rf{|xi|},
\bel{|X|u}\|X(\cd)\|_\sX\les K\|\xi(\cd)\|_\sX\les K\(\|\f(\cd)\|_\sX
+\|u(\cd)\|_\sU\).\ee
We call \rf{X=u} the variation of constants formula for the state $X(\cd)$.
From the above, we see that under (A1), for any control $u(\cd)\in\sU$, the state
equation \rf{state} is well-posed in $\sX$. Thus, the running cost in \rf{cost}
is well-defined. However, the terminal cost is still not necessarily defined. We
need the state process $X(\cd)$ to be continuous at $t=T$. To achieve this, we need
a little more assumption which we now introduce.

\ms

{\bf(A2)} \rm Let (A1) hold, and in addition, there exists a modulus of continuity
$\o(\cd)$ and some $\d_0\in(0,T]$,
$$|A(T,s)-A(t,s)|+|B(T,s)-B(t,s)|+|\f(T)-\f(t)|\les\o(T-t),\qq t\in[T-\d_0,T],\q
(t,s)\in\D.$$

\ms

{\bf(A3)} \rm $\b>{1\over2}$.

\ms

We have the following result.

\bt{continuity T} \sl {\rm(i)} Let {\rm(A2)--(A3)} hold. Then for any $s\in[0,T)$, $t\mapsto\F(t,s)$ is continuous at $t=T$.

\ms

{\rm(ii)} Let {\rm(A2)--(A3)} hold. Then for any control $u(\cd)\in\sU$,
the corresponding state process $X(\cd)$ is continuous at $t=T$.

\et
\it Proof. \rm (i) By Corollary \ref{corollary 2.2} (ii), we can get (i).
\ms

(ii) Now, we let (A2)--(A3) hold. Then for any $u(\cd)\in\sU$,
$$\ba{ll}
\ns\ds|X(t)|\les|\psi(t)|+\int_0^t|\Psi(t,s)u(s)|ds\les|\f(t)|+\int_0^t|\F(t,s)\f(s)|ds
+\(\int_0^t|\Psi(t,s)|^2ds\)^{1\over2}\|u(\cd)\|_\sU\\
\ns\ds\les|\f(t)|+\(\int_0^t|\F(t,s)|^2ds\)^{1\over2}\|\f(\cd)\|_\sX+K\(\int_0^t{ds\over
(t-s)^{2(1-\b)}}\)^{1\over2}\|u(\cd)\|_\sU\\
\ns\ds\les|\f(t)|+K\(\int_0^t{ds\over(t-s)^{2(1-\b)}}\)^{1\over2}\|\f(\cd)\|_\sX+{Kt^{\b-\frac{1}{2}}\over\sqrt{2\b-1}}\|u(\cd)\|_\sU\\
\ns\ds\les|\f(t)|+{Kt^{\b-\frac{1}{2}}\over\sqrt{2\b-1}}\|\f(\cd)\|_\sX+{Kt^{\b-\frac{1}{2}}\over\sqrt{2\b-1}}\|u(\cd)\|_\sU\les|\f(t)|+K\(\|\f(\cd)\|_\sX+\|u(\cd)\|_\sU\).\ea$$
Thus, $X(t)$ is defined at all points where $\f(t)$ is defined. To obtain the continuity of $X(\cd)$ at $t=T$, since $\f(\cd)$ is continuous at $T$, it suffices to obtain the
continuity of the following expression at $t=T$:
$$\int_0^t{A(t,\t)X(\t)+B(t,\t)u(\t)\over(t-\t)^{1-\b}}d\t.$$
Since
$$\ba{ll}
\ns\ds|A(t,\t)X(\t)+B(t,\t)u(\t)|\les K\(|X(\t)|+|u(\t)|\),\qq(t,\t)\in\D,\\
\ns\ds |A(T,\t)X(\t)+B(T,\t)u(\t)-A(t,\t)X(\t)-B(t,\t)u(\t)|\les\o(T-t)\(|X(\t)|+|u(\t)|\),\\
\ns\ds\qq\qq\qq\qq\qq\qq\qq\qq\qq\qq t\in[T-\d_0,T]$, $(t,\t)\in\D,\ea$$
with $|X(\cd)|+|u(\cd)|\in L^2(0,T)$, by Corollary \ref{corollary 2.2} (ii), we obtain the
continuity. \endpf

\ms

Note that if $\b\les{1\over2}$, then in general, we do not expect to have a continuity of the state process at $t=T$ for some control $u(\cd)\in\sU$. The following example illustrates this.

\bex{} \rm Let $T=1$, $A(\cd,\cd)=0$, $B(\cd,\cd)=1$, and $\f(\cd)=0$. Then we have
$$X(t)=\int_0^t{u(s)\over(t-s)^{1-\b}}ds,\qq t\in[0,1].$$
Let
$$u(s)={{\bf1}_{[{1\over2},1)}(s)\over(1-s)^{1\over2}\log(1-s)},\qq s\in[0,1).$$
Then,
$$\int_0^1|u(s)|^2ds=\int_{1\over2}^1{ds\over(1-s)[\log(1-s)]^2}={1\over\log2}.$$
Thus, $u(\cd)\in\sU$. However,
$$X(1)=\int_0^1{u(s)\over(1-s)^{1-\b}}ds=\int_{1\over2}^1{ds\over(1-s)^{{3\over2}-\b}\log(1-s)}=\infty,\qq\forall\b\les{1\over2}.$$

\ex
Having the above result, we see that under (A2)--(A3), for any $u(\cd)\in\sU$, the cost functional is well-defined, and therefore, Problem (P) is well-formulated.

\section{Open-loop Optimal Control.}

In this section, we will present the unique existence of open-loop optimal control and its characterizations. Let us first introduce the the following operators:
\bel{Th}\left\{\2n\ba{ll}
\ds(\Th u)(t)=\int_0^t\Psi(t,s)u(s)ds,\qq u\in\sU,~t\in[0,T],\\
\ns\ds\Th_Tu=\int_0^T\Psi(T,s)u(s)ds=(\Th u)(T),\qq u\in\sU.\ea\right.\ee
Then
\bel{X=}X(t)=\psi(t)+(\Th u)(t),\qq t\in[0,T];\qq\qq X(T)=\psi(T)+\Th_Tu.\ee
From \rf{|Psi|}, we see that (only need $\b\in(0,1)$) for any $u(\cd)\in\sU,$
\bel{Th u}\|\Th u\|_\sX=\(\int_0^T\Big|\int_0^t\Psi(t,s)u(s)ds\Big|^2dt\)^{1\over2}
\les K\[\int_0^T\(\int_0^t{|u(s)|\over(t-s)^{1-\b}}ds\)^2dt\]^{1\over2}\les
K\|u(\cd)\|_\sU,\ee
and (noting $\b\in({1\over2},1)$)
\bel{Th_T}|\Th_Tu|=\Big|\int_0^T\Psi(T,s)u(s)ds\Big|\les K\int_0^T{|u(s)|\over(T-s)^{1-\b}}ds\les K\|u(\cd)\|_\sU.\ee
Thus, $\Th\in\sL(\sU;\sX)$ and $\Th_T\in\sL(\sU;\dbR^n)$. Consequently, their adjoint operators $\Th^*\in\sL(\sX;\sU)$ and $\Th_T^*\in\sL(\dbR^n;\sU)$ are well-defined. Let us identify them as follows. For any $X(\cd)\in\sX$,
$$\ba{ll}
\ns\ds\lan X(\cd),(\Th u)(\cd)\ran_\sX=\int_0^T\lan X(t),(\Th u)(t)\ran dt=\int_0^T\lan X(t),\int_0^t\Psi(t,s)u(s)ds\ran dt\\
\ns\ds=\int_0^T\int_0^t\lan X(t),\Psi(t,s)u(s)\ran dsdt=\int_0^T\int_s^T\lan X(t),\Psi(t,s)u(s)\ran dtds\\
\ns\ds=\int_0^T\lan\int_s^T\Psi(t,s)^\top X(t)dt,u(s)\ran ds=\lan\int_s^T\Psi(t,s)^\top X(t)dt,u(s)\ran_\sU=\lan(\Th^*X)(\cd),u(\cd)\ran_\sU.\ea$$
This gives
\bel{Th*}(\Th^*X)(s)=\int_s^T\Psi(t,s)^\top X(t)dt,\qq X(\cd)\in\sX.\ee
From Corollary \ref{corollary 2.2}, (i), we see that (only need $\b\in(0,1)$)
$$\ba{ll}
\ns\ds\(\int_0^T\Big|\int_s^T\Psi(t,s)^\top X(t)dt\Big|^2ds\)^{1\over2}
\les K\[\int_0^T\(\int_s^T{|X(t)|\over(t-s)^{1-\b}}dt\)^2ds\]^{1\over2}\les K\|X(\cd)\|_\sX.\ea$$
Likewise (noting $\b\in({1\over2},1)$),
\bel{Th_T*}(\Th_T^*x)(s)=\Psi(T,s)^\top x,\qq\forall x\in\dbR^n,\ee
with
$$\(\int_0^T|(\Th^*_Tx)(s)|^2ds\)^{1\over2}\les K|x|\(\int_0^T{ds\over(T-s)^{2(1-\b)}}
\)^{1\over2}=K\({T^{2\b-1}\over2\b-1}\)^{1\over2}|x|.$$

\br{UU} \rm By \rf{|Psi|}, noting $\b>{1\over2}$, we see that for $X(\cd)\in\sX$,
$$|(\Th^*X)(s)|\les\int_s^T|\Psi(t,s)||X(t)|dt\les K\int_s^T{|X(t)|\over
(t-s)^{1-\b}}dt\les K{(T-s)^{\b-{1\over2}}\over(2\b-1)^{1\over2}}\|X(\cd)\|_\sX,\q s\in[0,T],$$
which is bounded. But,
$$|(\Th^*_Tx)(s)|\les{K|x|\over(T-s)^{1-\b}},\qq s\in[0,T),$$
which may be unbounded. This is different from the case that \rf{K<infty} holds as in
\cite{Prichard-You 1996}. Consequently, we need to modify the technique used there
so that it works for us.

\er
Now, we would like to represent the cost functional. We will use $\BQ$ to denote the bounded operator from $\sX$ to itself induced by $Q(\cd)$, and so on. Then
\bel{rep}\ba{ll}
\ns\ds J(u(\cd))=\lan\BQ X,X\ran+2\lan\BS X,u\ran+\lan\BR u,u\ran+2\lan\Bq,X\ran+2\lan\BBrho,u\ran+\lan\BG X(T),X(T)\ran+2\lan\Bg,X(T)\ran\\
\ns\ds\qq\q=\lan\BQ(\psi+\Th u),\psi+\Th u\ran+2\lan\BS(\psi+\Th u),
u\ran+\lan\BR u,u\ran+2\lan\Bq,\psi+\Th u\ran+2\lan\BBrho,u\ran\\
\ns\ds\qq\qq+\lan\BG[\psi(T)+\Th_Tu],\psi(T)+\Th_Tu\ran+2\lan\Bg,\psi(T)+\Th_Tu\ran\\
\ns\ds\qq\q=\lan\BBL u,u\ran+2\lan\BBl_1,u\ran+\l_0,\ea\ee
where
$$\left\{\2n\ba{ll}
\ds\BBL=\Th^*\BQ\Th+\BS\Th+\Th^*\BS^*+\BR+\Th_T^*\BG
\Th_T,\\
\ns\ds\BBl_1=\Th^*\BQ\psi+\BS\psi+\Th^*\Bq+\BBrho+
\Th_T^*\BG\psi(T)+\Th_T^*\Bg,\\
\ns\ds\l_0=\lan\BQ\psi,\psi\ran+2\lan\Bq,\psi\ran+
\lan\BG\psi(T),\psi(T)\ran+2\lan\Bg,\psi(T)\ran.\ea\right.$$
The above shows that our Problem (P) can be regarded as an optimization of the functional $J(u(\cd))$ on the Hilbert space $\sU$. In order the functional $u(\cd)\mapsto J(u(\cd))$ to be bounded from below, it
is necessary that $\BBL\ges0$. In what follows, we want $J(u(\cd))$ to admit a
unique minimum. To guarantee this, we may assume the following stronger
condition:
\bel{L_2>0*}\BBL\ges\d,\ee
for some $\d>0$. Since it is not the main theme of this paper to discuss the
conditions under which \rf{L_2>0*} holds, we are satisfied to assume
proper conditions to guarantee \rf{L_2>0*}. For simplicity, we introduce the
following standard assumption:

\ms

{\bf(A4)} Let
$$\left\{\2n\ba{ll}
\ns\ds Q(\cd)\in L^\infty(0,T;\dbS^n),\qq S(\cd)\in L^\infty(0,T;\dbR^{m\times n}),\qq
R(\cd)\in L^\infty(0,T;\dbS^m),\\
\ns\ds q(\cd)\in L^2(0,T;\dbR^n),\qq\rho(\cd)\in L^2(0,T;\dbR^m),\qq
 G\in\dbS^n,\qq g\in\dbR^n,\ea\right.$$
and the following holds:
\bel{classic}R(t)\ges\d,\q Q(t)-S(t)^\top R(t)^{-1}S(t)
\ges0,\q G\ges0,\q t\in[0,T],\ee
for some $\d>0$.

\ms

Now, we have the following result for Problem (P).

\bt{tho1} \sl Suppose that {\rm(A2)--(A4)} hold. Then Problem {\rm(P)} admits a unique
open-loop optimal pair $(\bar X(\cd),\bar u(\cd))\in\sX\times\sU$. Moreover, the following
relation is satisfied:
\bel{bar u1}\ba{ll}
\ns\ds\bar u(t)=-R(t)^{-1}\Big\{\Psi(T,t)^\top\(G\bar X(T)+g\)+S(t)\bar X(t)+\rho(t)\\
\ns\ds\qq\qq\qq+\int_t^T
\Psi(s,t)^\top\(Q(s)\bar X(s)+S(s)^\top\bar u(s)+q(s)\)ds\Big\},\qq\ae~t\in[0,T].\ea\ee

\et

\it Proof. \rm Under (A4), one has \rf{L_2>0*}. Thus, from the representation \rf{rep} of the cost functional, we see that $u(\cd)\mapsto J(u(\cd))$ admits a unique minimum
$\bar u(\cd)$ which is given by the solution to the following:
$$\ba{ll}
\ns\ds0=\BBL\bar u+\BBl_1=\(\Th^*\BQ\Th+\BS\Th+\Th^*\BS^*+\BR+\Th_T^*\BG
\Th_T\)\bar u+\Th^*\BQ\psi+\BS\psi+\Th^*\Bq+\BBrho+
\Th_T^*\BG\psi(T)+\Th_T^*\Bg\\
\ns\ds\q=\Th^*\BQ(\psi+\Th\bar u)+\BS(\psi+\Th\bar u)+\Th^*\BS^*\bar u+\BR\bar u+\Th_T^*\BG(\psi(T)+\Th_T\bar u)+\Th^*\Bq+\BBrho+\Th_T^*\Bg\\
\ns\ds\q=\Th_T^*[\BG\bar X(T)+\Bg]+\Th^*(\BQ\bar X+\BS^*\bar u+\Bq)+\BR\bar u+\BS\bar X+\BBrho.\ea$$
Thus,
\bel{bar u0}\bar u=-\BR^{-1}\(\Th_T^*[\BG\bar X(T)+\Bg]+\BS\bar X+\BBrho+\Th^*(\BQ\bar X+\BS^*\bar u+\Bq)\),\ee
which is the same as \rf{bar u1}. \endpf

\ms

Let us take a closer look at the above result. Note that
$$\Psi(s,t)^\top={B(s,t)^\top\over(s-t)^{1-\b}}
+\int_t^s{B(\t,t)^\top\F(s,\t)^\top\over(\t-t)^{1-\b}}d\t,\qq(s,t)\in\D.$$
Thus,
$$\ba{ll}
\ns\ds\q\Psi(T,t)^\top\(G\bar X(T)+g\)+\int_t^T
\Psi(s,t)^\top\[Q(s)\bar X(s)+S(s)^\top\bar u(s)+q(s)\]ds\\
\ns\ds={B(T,t)^\top[G\bar X(T)+g]\over(T-t)^{1-\b}}
+\int_t^T{B(s,t)^\top\F(T,s)^\top[G\bar X(T)+g]\over(s-t)^{1-\b}}ds\\
\ns\ds\q+\int_t^T{B(s,t)^\top\over(s-t)^{1-\b}}\Big\{Q(s)\bar X(s)+S(s)^\top\bar u(s)+q(s)+\int_s^T\F(\t,s)^\top\big[Q(\t)\bar X(\t)+S(\t)^\top\bar u(\t)+q(\t)\big]d\t\Big\}ds\ea$$
$$\ba{ll}
\ns\ds={B(T,t)^\top[G\bar X(T)+g]\over(T-t)^{1-\b}}+\int_t^T{B(s,t)^\top\over(s-t)^{1-\b}}\Big\{\F(T,s)^\top[G\bar X(T)+g]+Q(s)\bar X(s)+S(s)^\top\bar u(s)+q(s)\\
\ns\ds\qq\qq\qq\qq\qq\qq\qq\qq+\int_s^T\F(\t,s)^\top\big[Q(\t)\bar X(\t)+S(\t)^\top\bar u(\t)+q(\t)\big]d\t\Big\}ds,\q t\in[0,T].\ea$$
Consequently, we have the following relation for the optimal control $\bar u(\cd)$:
\bel{bar u2}\ba{ll}
\ns\ds\bar u(t)=-R(t)^{-1}\[{B(T,t)^\top[G\bar X(T)+g]\over(T-t)^{1-\b}}+S(t)\bar X(t)+\rho(t)\\
\ns\ds\qq\qq+\int_t^T{B(s,t)^\top\over(s-t)^{1-\b}}\(\F(T,s)^\top[G\bar X(T)+g]+Q(s)\bar X(s)+S(s)^\top\bar u(s)+q(s)\\
\ns\ds\qq\qq\qq\qq\qq+\int_s^T\F(\t,s)^\top\big[Q(\t)\bar X(\t)+S(\t)^\top\bar u(\t)+q(\t)\big]d\t\)ds\],\qq\ae~t\in[0,T],\ea\ee
with $\bar X(\cd)$ being the optimal state trajectory.

\ms

On the other hand, we know that optimal control can also be characterized by the variational method. The following is the result for Problem (P) from a different angle.

\bt{maximum principle} \sl Let {\rm(A2)--(A4)} hold. Then Problem {\rm(P)} admits a unique open-loop optimal pair $(\bar X(\cd),\bar u(\cd))$ such that
\bel{bar u3}\bar u(t)=-R(t)^{-1}\[\int_t^T{B(s,t)^\top Y(s)\over(s-t)^{1-\b}}ds+S(t)\bar X(t)+\rho(t)+{B(T,t)^\top[G\bar X(T)+g]\over(T-t)^{1-\b}}\],\qq\ae~t\in[0,T],\ee
where $Y(\cd)$ is the solution to the following adjoint equation:
\bel{Y}Y(t)=Q(t)\bar X(t)+S(t)^\top\bar u(t)+q(t)+{A(T,t)^\top[G\bar X(T)+g]\over(T-t)^{1-\b}}+\int_t^T{A(s,t)^\top Y(s)\over(s-t)^{1-\b}}ds,
\qq\ae~t\in[0,T].\ee

\et

\it Proof. \rm Let $(\bar X(\cd),\bar u(\cd))$ be the optimal pair. Then for any $u(\cd)\in\sU$, we have
$$\ba{ll}
\ns\ds 0=\lim_{\e\to0}{1\over2\e}\[J(\bar u(\cd)+\e u(\cd))-J(\bar u(\cd))\]\\
\ns\ds=\int_0^T\(\lan Q(t)\bar X(t)+S(t)^\top\bar u(t)+q(t),X(t)\ran+\lan S(t)\bar X(t)+R(t)\bar u(t)+\rho(t),u(t)\ran\)dt+\lan G\bar X(T)+g,X(T)\ran,\ea$$
where
$$X(t)=\int_0^t{A(t,s)X(s)+B(t,s)u(s)\over(t-s)^{1-\b}}ds,\qq t\in[0,T].$$
Thus,
$$\ba{ll}
\ns\ds\lan G\bar X(T)+g,X(T)\ran=\int_0^T\(\lan{A(T,t)^\top[G\bar X(T)+g]\over(T-t)^{1-\b}},X(t)\ran +\lan{B(T,t)^\top[G\bar X(T)+g]\over(T-t)^{1-\b}},u(t)\ran\)dt.\ea$$
This yields
$$\ba{ll}
\ns\ds0=\int_0^T\(\lan Q(t)\bar X(t)+S(t)^\top\bar u(t)+q(t)+{A(T,t)^\top[G\bar X(T)+g]\over(T-t)^{1-\b}},X(t)\ran\\
\ns\ds\qq\qq+\lan S(t)\bar X(t)+R(t)\bar u(t)+\rho(t)+{B(T,t)^\top[G\bar X(T)+g]\over(T-t)^{1-\b}},u(t)\ran\)dt.\ea$$
Now, we let $Y(\cd)$ be the solution to the following:
$$Y(t)=\gamma(t)+\int_t^T{A(s,t)^\top Y(s)\over(s-t)^{1-\b}}ds,\qq\ae~t\in[0,T],$$
with
$$\gamma(t)=Q(t)\bar X(t)+S(t)^\top\bar u(t)+q(t)+{A(T,t)^\top[G\bar X(T)+g]\over(T-t)^{1-\b}},\qq\ae~t\in[0,T].$$
Then for any $u\in\sU,$
$$\ba{ll}
\ns\ds0=\int_0^T\(\lan Y(t)-\int_t^T{A(s,t)^\top Y(s)\over(s-t)^{1-\b}}ds,X(t)\ran\\
\ns\ds\qq\qq+\lan S(t)\bar X(t)+R(t)\bar u(t)+\rho(t)+{B(T,t)^\top[G\bar X(T)+g]\over(T-t)^{1-\b}},u(t)\ran\)dt\\
\ns\ds\q=\int_0^T\lan Y(t),X(t)\ran dt-\int_0^T\lan Y(s),\int_0^s{A(s,t)X(t)\over(s-t)^{1-\b}}dt\ran ds\\
\ns\ds\qq\qq+\int_0^T\lan S(t)\bar X(t)+R(t)\bar u(t)+\rho(t)+{B(T,t)^\top[G\bar X(T)+g]\over(T-t)^{1-\b}},u(t)\ran dt\\
\ns\ds\q=\int_0^T\(\lan Y(t),\int_0^t{B(t,s)u(s)\over(t-s)^{1-\b}}ds\ran+\lan S(t)\bar X(t)+R(t)\bar u(t)+\rho(t)+{B(T,t)^\top[G\bar X(T)+g]\over(T-t)^{1-\b}},u(t)\ran\)dt\\
\ns\ds\q=\int_0^T\lan\int_t^T{B(s,t)^\top Y(s)\over(s-t)^{1-\b}}ds+S(t)\bar X(t)+R(t)\bar u(t)+\rho(t)+{B(T,t)^\top[G\bar X(T)+g]\over(T-t)^{1-\b}},u(t)\ran dt.\ea$$
Hence,
$$\int_t^T{B(s,t)^\top Y(s)\over(s-t)^{1-\b}}ds+S(t)\bar X(t)+R(t)\bar u(t)+\rho(t)+{B(T,t)^\top[G\bar X(T)+g]\over(T-t)^{1-\b}}=0,\qq\ae~t\in[0,T].$$
Then \rf{bar u3} follows. \endpf

\ms

The above is actually the Pontryagin type maximum principle. Comparing \rf{bar u2} and \rf{bar u3}, we see that they coincide if the following is true:
\bel{Y*}\ba{ll}
\ns\ds Y(s)=\F(T,s)^\top[G\bar X(T)+g]+Q(s)\bar X(s)+S(s)^\top\bar u(s)+q(s)\\
\ns\ds\qq\qq\qq\qq\qq+\int_s^T\F(\t,s)^\top\big[Q(\t)\bar X(\t)+S(\t)^\top\bar u(\t)+q(\t)\big]d\t,\qq\ae~s\in[0,T].\ea\ee
This can be shown as follows. By \rf{Phi*}, we have
\bel{Phi T}\F(t,s)^\top={A(t,s)^\top\over(t-s)^{1-\b}}
+\int_s^t{A(\t,s)^\top\F(t,\t)^\top\over(\t-s)^{1-\b}}d\t,\qq \ae~t\in[s,T].\ee
Denote
$$z(t)=Q(t)\bar X(t)+S(t)^\top\bar u(t)+q(t),\q\ae~t\in[0,T];\qq\z=G\bar X(T)+g.$$
Then, we have
$$Y(t)=z(t)+{A(T,t)^\top\z\over(T-t)^{1-\b}}
+\int_t^T{A(s,t)^\top Y(s)\over(s-t)^{1-\b}}ds,\qq\ae~t\in[0,T],$$
and we need to check:
$$Y(t)=z(t)+\F(T,t)^\top\z+\int_t^T\F(s,t)^\top z(s)ds,\qq\ae~t\in[0,T].$$
This can be checked as follows:
$$\ba{ll}
\ns\ds\int_t^T{A(s,t)^\top Y(s)\over(s-t)^{1-\b}}ds=\int_t^T{A(s,t)^\top
\big[z(s)+\F(T,s)^\top\z+\int_s^T\F(\t,s)^\top z(\t)d\t\big]\over
(s-t)^{1-\b}}ds\\
\ns\ds=\int_t^T{A(s,t)^\top z(s)\over(s-t)^{1-\b}}ds+\int_t^T{A(s,t)^\top\F(T,s)^\top\z\over(s-t)^{1-\b}}ds
+\int_t^T\(\int_t^s{A(\t,t)^\top\F(s,\t)^\top\over(\t-t)^{1-\b}}d\t\)z(s)ds\\
\ns\ds=\int_t^T\F(s,t)^\top z(s)ds+\int_t^T{A(s,t)^\top\F(T,s)^\top\z\over(s-t)^{1-\b}}ds\\
\ns\ds=Y(t)-z(t)-\(\F(T,t)^\top-\int_t^T{A(s,t)^\top\F(T,s)^\top\over(s-t)^{1-\b}}ds\)
\z=Y(t)-z(t)-{A(T,t)^\top\z\over(T-t)^{1-\b}},\qq\ae~t\in[0,T].\ea$$
Hence, \rf{bar u2} and \rf{bar u3} are equivalent.

\ms

\section{Causal State Feedback Representation.}

In the relation \rf{bar u2} (or \rf{bar u3}) for the open-loop optimal control, the current-time value $\bar u(t)$ of the optimal control $\bar u(\cd)$ is given in terms of the future-time values $\{\bar X(s)\bigm|s\in[t,T]\}$ of the corresponding optimal state trajectory $\bar{X}(\cd)$. Practically, this is not realizable. Thus, our next goal is to seek a {\it causal state feedback representation} of optimal control, by which we mean that the value $\bar u(t)$ of $\bar u(\cd)$ at time $t$ can be written in terms of $\{\bar X(s)\bigm|s\in[0,t]\}$ and $\bar X^a(t)$ for some non-anticipating auxiliary process $\bar X^a(\cd)$. Recall that for standard LQ problems of differential equations, optimal control could admit a closed-loop representation by means of differential Riccati equations. Here, we borrow the idea from \cite{Prichard-You 1996}, but with more straightforward approach which more naturally reveals the essence of the problem.

\ms

Since our cost functional contains the cross term of the state and control, as well as a linear term in the control, we would like to make a reduction first.
Let
\bel{u=v}u(t)=v(t)-R(t)^{-1}\[S(t)X(t)+\rho(t)\],\qq\ae~ t\in[0,T].\ee
Then, the state equation becomes
\bel{state-new}X(t)=\h\f(t)+\int_0^t{\h A(t,s)X(s)+B(t,s)v(s)\over(t-s)^{1-\b}}ds,\qq\ae~ t\in[0,T],\ee
where
$$\h A(t,s)=A(t,s)-B(t,s)R(s)^{-1}S(s),\qq\h\f(t)=\f(t)-\int_0^t{B(t,s)R(s)^{-1}\rho(s)
\over(t-s)^{1-\b}}ds,\qq\ae~(t,s)\in\Delta.$$
Also, the running cost rate becomes
$$\ba{ll}
\ns\ds\q\lan QX,X\ran+2\lan SX,u\ran+\lan Ru,u\ran+2\lan q,X\ran+2\lan\rho,u\ran\\
\ns\ds=\lan(Q-S^\top R^{-1}S)X,X\ran+\lan Rv,v\ran+2\lan q-S^\top R^{-1}\rho,X\ran-\lan R^{-1}\rho,\rho\ran\\
\ns\ds\equiv\lan\h QX,X\ran+\lan Rv,v\ran+2\lan\h q,X\ran-\lan R^{-1}\rho,\rho\ran.\ea$$
Thus, if we define
$$\h J(v(\cd))=\int_0^T\(\lan\h Q(t)X(t),X(t)\ran+\lan R(t)v(t),v(t)\ran+2\lan\h q(t),X(t)\ran\)dt+\lan GX(T),X(T)\ran+2\lan g,X(T)\ran,$$
then the corresponding optimal control problem is equivalent to the original one.
Namely, $(\bar X(\cd),\bar u(\cd))$ is the open-loop optimal pair of the original
problem if and only if $(\bar X(\cd),\bar v(\cd))$ is the optimal pair of the
reduced problem with $\bar u(\cd)$ and $\bar v(\cd)$ being related by the following:
\bel{u=v}\bar u(t)=\bar{v}(t)-R(t)^{-1}\big[S(t)\bar X(t)+\rho(t)\big],\qq\ae~t\in[0,T].\ee
Therefore, if the optimal control $\bar v(\cd)$ of the reduced problem has a causal state feedback representation, then so is $\bar u(\cd)$. Hence, for simplicity, we
introduce the following hypothesis.

\ms

{\bf(A5)} Let (A2)--(A4) hold with
\bel{S=0}S(t)=0,\q\rho(t)=0,\qq t\in[0,T].\ee

Under (A5), we have
\bel{L*}\left\{\2n\ba{ll}
\ds\BBL=\Th^*\BQ\Th+\Th_T^*\BG\Th_T+\BR,\\
\ns\ds\BBl_1=\Th^*\BQ\psi+\Th^*\Bq+\Th_T^*\BG\psi(T)+\Th_T^*\Bg,\\
\ns\ds\l_0=\lan\BQ\psi,\psi\ran+2\lan\Bq,\psi\ran+
\lan\BG\psi(T),\psi(T)\ran+2\lan\Bg,\psi(T)\ran,\ea\right.\ee
and the open-loop optimal control is given by
\bel{bar u0}\bar u=-\BR^{-1}\(\Th^*(\BQ\bar X+\Bq)+\Th_T^*[\BG\bar X(T)+\Bg]\)
=-\BR^{-1}\big[\Th^*\BQ\bar X+\Th^*_T\BG\bar X(T)\big]-\BR^{-1}\big[\Th^*\Bq+\Th^*_T\Bg\big],\ee
or more precisely,
\bel{bar u(t)}\ba{ll}
\ds\bar u(t)=-R(t)^{-1}\[\int_t^T\Psi(\t,t)^\top Q(\t)\bar X(\t)d\t+\Psi(T,t)^\top G\bar X(T)\]\\
\ds\qq\q-R(t)^{-1}\[\int_t^T\Psi(\t,t)^\top q(\t)d\t+\Psi(T,t)^\top g\],\qq\ae~t\in[0,T].\ea\ee
In the above, $\Psi(\cd\,,\cd)$, which is defined by \rf{psi}, characterizes the control system and $Q(\cd)$, $q(\cd)$, $G$ and $g$ are all known a priori. The only unrealistic  terms are $\{\bar X(\t)\bigm|\t\in[t,T]\}$ and $\bar X(T)$ on the right-hand side of the above, since at the time $t$ of determining the value $\bar u(t)$ of $\bar u(\cd)$, these are not available.

\ms

The idea of getting a feasible representation of the optimal control is to introduce the following simple decomposition for the state trajectory: For any $\si\in[0,T)$,
$$\ba{ll}
\ns\ds\bar X(t)=\psi(t)+\int_0^t\Psi(t,s)\bar u(s){\bf1}_{[0,\si)}(s)ds+\int_0^t\Psi(t,s)\bar u(s){\bf1}_{[\si,T]}(s)ds\\
\ns\ds\qq=\psi(t)+\int_0^{t\land\si}\Psi(t,s)\bar u(s)ds+\int_{t\land\si}^t\Psi(t,s)\bar u(s)ds\equiv\bar X_\si(t)+\int_{t\land\si}^t\Psi(t,s)\bar u(s)ds,\qq\ae~t\in[0,T].\ea$$
Also,
$$\bar X(T)=\psi(T)+\int_0^\si\Psi(T,s)\bar u(s)ds+\int_\si^T\Psi(T,s)\bar u(s)ds\equiv\bar X^a(\si)+\int_\si^T\Psi(T,s)\bar u(s)ds.$$
Here, $\bar X_\si(t)$ and $\bar X^a(\si)$ do not use the information of $\bar u(\cd)$ beyond $\si$. With such a decomposition, we can rewrite \rf{bar u(t)} as follows:
\bel{bar u(t)*}\ba{ll}
\ds\bar u(t)=-R(t)^{-1}\[\int_t^T\Psi(\t,t)^\top Q(\t)\(\bar X_\si(\t)+\int_{\t\land\si}^\t\Psi(\t,s)\bar u(s)ds\)d\t\\
\ds\qq\qq\qq\qq+\Psi(T,t)^\top G\(\bar X^a(\si)+\int_\si^T\Psi(T,s)\bar u(s)ds\)\]\\
\ds\qq\q-R(t)^{-1}\[\int_t^T\Psi(\t,t)^\top q(\t)d\t+\Psi(T,t)^\top g\],\qq\ae~t\in[0,T].\ea\ee
By letting $\si=t$ in the above, we obtain
\bel{bar u(t)**}\ba{ll}
\ns\ds\bar u(t)=-R(t)^{-1}\[\int_t^T\Psi(\t,t)^\top Q(\t)\bar X_t(\t)d\t+\Psi(T,t)^\top G\bar X^a(t)\]\\
\ns\ds\qq\q-R(t)^{-1}\[\int_t^T\Psi(\t,t)^\top Q(\t)\int_{\t\land t}^\t\Psi(\t,s)\bar u(s)ds d\t+\Psi(T,t)^\top G\int_t^T\Psi(T,s)\bar u(s)ds\]\\
\ns\ds\qq\q-R(t)^{-1}\[\int_t^T\Psi(\t,t)^\top q(\t)d\t+\Psi(T,t)^\top g\],\qq\ae~t\in[0,T].\ea\ee
According to the definition of $(\bar X_t(\cd),\bar X^a(t))$, no information of $\bar u(\cd)$ beyond $t$ is used, or they are non-anticipating. Therefore, our goal is to rewrite the last two terms on the right-hand side to be non-anticipating.

\ms

To achieve our goal, let us introduce a family of projection operators $\Pi_\si:\sU\to\sU$ by the following:
\bel{pi}[\Pi_\si u(\cd)](t)={\bf1}_{[0,\si)}(t)u(t),\qq t\in[0,T],\ee
with $\si\in[0,T)$ being the parameter. Clearly,
\bel{I-pi}[(I-\Pi_\si)u(\cd)](t)={\bf1}_{[\si,T]}(t)u(t),\qq t\in[0,T].\ee
Both $\Pi_\si$ and $I-\Pi_\si$ are idempotents on $\sU$. Moreover, for any
$M(\cd)\in L^\infty(0,T;\dbR^{k\times m})$,
$$M(t)[\Pi_\si u(\cd)](t)=M(t){\bf1}_{[0,\si)}(t)u(t)={\bf1}_{[0,\si)}(t)M(t)u(t)=
[\Pi_\si M(\cd)u(\cd)](t),\qq t\in[0,T].$$
This means that $\Pi_\si$ commutes with multiplication operators. We shall call $\Pi_\si$  a {\it causal projection}. Now, let
$$\sU_\si=\sR\big(I-\Pi_\si\big)=(I-\Pi_\si)\sU,$$
which is a closed subspace of $\sU,$ and define a parameterized operator
$\BBL_\si\in\sL(\sU_\si)$ by
\bel{L_si}\BBL_\si=(I-\Pi_\si)\BBL\big|_{\sU_\si}.\ee
Clearly, $(I-\Pi_\si)\BBL(I-\Pi_\si)$ is a natural extension of $\BBL_\si$, with the value being 0 on $\sU_\si^\perp$.
%
%Indeed, let $u_n\in\sU^{+}_\xi$ with $u_n\rightarrow u,$ as $n\rightarrow\infty.$ Then there exists $v_n\in\sU$ such that $u_n=(I-\pi_{\xi})v_n.$ Then,
%
%$$u\leftarrow(I-\pi_{\xi})v_n =(I-\pi_{\xi})(I-\pi_{\xi})v_n\rightarrow(I-\pi_{\xi})u.$$
%
%Thus, $(I-\pi_{\xi})u=u,$ $u\in\sU.$ Then, $u\in\sU^{+}_\xi.$ Thus $\sU^{+}_\xi$ is a closed subspace of $\sU.$
%
We have the following simple lemma.

\bl{lemma operator} \sl Suppose that {\rm(A5)} holds. Then for any given $\si\in[0,T)$, the operator $\BBL_\si\in\sL(\sU_\si)$ is self-adjoint and positive definite on $\sU_\si$. Moreover
\bel{<1/d}\BBL_\si^{-1}=(I-\Pi_\si)\BBL^{-1}\big|_{\sU_\si};\qq\|(\BBL_\si)^{-1}\|_{\sL(\sU_\si)}\les{1\over\d},\qq\forall\si\in[0,T).\ee

\el

\it Proof. \rm For any $u,v\in\sU_\si$, we have
$$\ba{ll}
\ns\ds \lan\BBL_\si u,v\ran_{\sU_\si}=\lan\BBL_\si u,v\ran_\sU=
\lan(I-\Pi_\si)\BBL u,v\ran_\sU=\lan\BBL u,(I-\Pi_\si)v\ran_\sU=\lan\BBL u,v\ran_\sU\\
\ns\ds\q\qq\qq=\lan u,\BBL v\ran_\sU=\lan(I-\Pi_\si)u,\BBL v\ran_\sU=
\lan u,\BBL_\si v\ran_\sU=\lan u,\BBL_\si v\ran_{\sU_\si}.\ea$$
Thus, $\BBL_\si$ is self-adjoint on $\sU_\si$. Moreover, there exists a constant $\d>0$ such that for any  $u\in\sU_\si$,
$$\lan\BBL_\si u,u\ran_{\sU_\si}=\lan\BBL(I-\Pi_\si)u,(I-\Pi_\si)u\ran_\sU
\ges\d\|(I-\Pi_\si)u\|^2_\sU=\d\|u\|^2_{\sU_\si}.$$
Consequently, $\BBL_\si$ is boundedly invertible. Applying the above to $(\BBL_\si)^{-1}u$, one has
$$\d\|(\BBL_\si)^{-1}u\|^2_{\sU_\si}\les\lan\BBL_\si(\BBL_\si)^{-1}u,(\BBL_\si)^{-1}u
\ran_{\sU_\si}=\lan u,(\BBL_\si)^{-1}u\ran_{\sU_\si}\les\|u\|_{\sU_\si}\|
(\BBL_\si)^{-1}u\|_{\sU_\si}.$$
Then the second estimate in \rf{<1/d} follows. \endpf

\ms

From the above, we see that
\bel{LL}\ba{ll}
\ds\BBL_\si(I-\Pi_\si)=(I-\Pi_\si)\BBL(I-\Pi_\si),\qq\BBL_\si^{-1}(I-\Pi_\si)
=(I-\Pi_\si)\BBL^{-1}(I-\Pi_\si),\\
\ns\ds(I-\Pi_\si)\BBL(I-\Pi_\si)\BBL^{-1}(I-\Pi_\si)=I-\Pi_\si,\q(I-\Pi_\si)
\BBL^{-1}(I-\Pi_\si)\BBL(I-\Pi_\si)=I-\Pi_\si.\ea\ee
Now, we introduce the following {\it auxiliary trajectory}:
\bel{X^a}X^a(t)=\psi(T)+\int_0^t\Psi(T,s)u(s)ds,\qq t\in[0,T],\ee
which catches the anticipating information of the free term $\psi(T)$ and the
dynamic system represented by $\Psi(T,\cd)$ which are assumed to be a priori known\footnote{In time-varying LQ problems, the differential Riccati equation is solved on the whole time interval and all information of the system and the cost functional through the coefficients and the weights are allowed to be used. The situation here is similar.}.
Note that no anticipating information of the control is involved. At the same time,
for each $\si\in[0,T)$, we introduce the following {\it Causal trajectory}:
\bel{X_si}X_\si(t)=\psi(t)+\int_0^t\Psi(t,s)(\Pi_\si u)(s)ds=\psi(t)+\int_0^{t\land\si}\Psi(t,s)u(s)ds,\qq t\in[0,T].\ee
This trajectory truncates the control up to time moment $\si$ which is allowed to
be smaller than $t$. It is clear that
\bel{X_si(t)}X_\si(t)=\left\{\2n\ba{ll}
\ds X(t),\qq t\in[0,\si],\\
\ns\ds\psi(t)+\int_0^\si\Psi(t,s)u(s)ds,\q t\in[\si,T];\ea\right.\ee
and
\bel{X^a(t)}X_t(T)=X^a(t),\qq t\in[0,T].\ee
We point out that both $X_\si(\cd)$ and $X^a(\cd)$ can be running at the same time as the system is running. Hence, $t\mapsto(X_t(\cd),X^a(t))$ is non-anticipating. We now prove the main result of this section.

\bt{representation1} \sl Let {\rm(A5)} hold. Then the unique open-loop optimal control $\bar u(\cd)$ admits the following representation:
\bel{bar u0*}\ba{ll}
\ns\ds\bar u(t)=-R(t)^{-1}\[I-(\BBL-\BR)(I-\Pi_t)\BBL^{-1}(I-\Pi_t)\]\(\big[\Th^*Q\bar X_t(\cd)\big](t)
+\big[\Th^*_TG\bar X^a(t)\big](t)\)\\
\ns\ds\qq\q-R(t)^{-1}\[I-(\BBL-\BR)(I-\Pi_t)\BBL^{-1}(I-\Pi_t)\]\(\big[\Th^*q(\cd)\big](t)
+\big[\Th^*_Tg\big](t)\),\qq\ae~ t\in[0,T],\ea\ee
where $\bar X_t(\cd)$ and $\bar X^a(\cd)$ are the causal and auxiliary trajectories corresponding to the optimal pair $(\bar X(\cd),\bar u(\cd))$.
\et

\it Proof. \rm By the definition of $X_\si(\cd)$ and $X^a(\cd)$, we have the following: for any $\si\in[0,T)$,
\bel{X=X_th}\ba{ll}
\ns\ds X(t)=\psi(t)+\int_0^t\Psi(t,s)\big\{(\Pi_\si u)(s)+[(I-\Pi_\si)u](s)\big\}ds=X_\si(t)+[\Th(I-\Pi_\si)u](t),\q t\in[0,T],\\
\ns\ds X(T)=\psi(T)+\int_0^T\Psi(T,s)\big\{(\Pi_\si u)(s)+[(I-\Pi_\si)u](s)\big\}ds=X^a(\si)+\Th_T(I-\Pi_\si)u.\ea\ee
Thus, the open-loop optimal control can be written as follows (see \rf{bar u0}):
\bel{bar u4}\ba{ll}
\ns\ds\bar u=-\BR^{-1}\big[\Th^*\BQ\bar X+\Th^*_T\BG\bar X(T)\big]-\BR^{-1}(\Th^*\Bq+\Th^*_T\Bg)\\
\ns\ds\q=-\BR^{-1}(\Th^*\BQ,\Th^*_T\BG)\begin{pmatrix}\bar X_\si+\Th(I-\Pi_\si)
\bar u\\ \bar X^a(\si)+\Th_T(I-\Pi_\si)\bar{u}\end{pmatrix}-\BR^{-1}(\Th^*\Bq+\Th^*_T\Bg)\\
\ns\ds\q=-\BR^{-1}(\Th^*\BQ,\Th^*_T\BG)\begin{pmatrix}\bar X_\si\\ \bar X^a(\si)\end{pmatrix}-\BR^{-1}(\BBL-\BR)(I-\Pi_\si)\bar u-\BR^{-1}(\Th^*\Bq+\Th^*_T\Bg).
\ea\ee
This leads to
$$\big[I+\BR^{-1}(\BBL-\BR)(I-\Pi_\si)\big]\bar u=-\BR^{-1}(\Th^*\BQ,\Th^*_T\BG)\begin{pmatrix}\bar X_\si\\ X^a(\si)\end{pmatrix}-\BR^{-1}(\Th^*\Bq+\Th^*_T\Bg),$$
which is equivalent to
$$\big[\BBL(I-\Pi_\si)+\BR\Pi_\si\big]\bar u=-(\Th^*\BQ,\Th^*_T\BG)\begin{pmatrix}\bar X_\si\\ X^a(\si)\end{pmatrix}-(\Th^*\Bq+\Th^*_T\Bg).$$
Applying $(I-\Pi_\si)$ to the above gives
$$\BBL_\si(I-\Pi_\si)\bar u=(I-\Pi_\si)\BBL(I-\Pi_\si)\bar u=-(I-\Pi_\si)(\Th^*\BQ,\Th^*_T\BG)\begin{pmatrix}\bar X_\si\\ X^a(\si)\end{pmatrix}-(I-\Pi_\si)(\Th^*\Bq+\Th^*_T\Bg).$$
Thus,
$$(I-\Pi_\si)\bar u=-\BBL_\si^{-1}(I-\Pi_\si)(\Th^*\BQ,\Th^*_T\BG)\begin{pmatrix}\bar X_\si\\ X^a(\si)\end{pmatrix}-\BBL_\si^{-1}(I-\Pi_\si)(\Th^*\Bq+\Th^*_T\Bg).$$
Substituting the above into \rf{bar u4}, one obtains
\bel{bar u5}\ba{ll}
\ns\ds\bar u=
%-\BR^{-1}(\Th^*\BQ,\Th^*_T\BG)\begin{pmatrix}\bar X_\si\\ \bar X^a(\si)\end{pmatrix}-\BR^{-1}(\BBL-\BR)(I-\Pi_\si)\bar u-\BR^{-1}(\Th^*\Bq+\Th^*_T\Bg)\\
%
%\ns\ds\q=
-\BR^{-1}(\Th^*\BQ,\Th^*_T\BG)\begin{pmatrix}\bar X_\si\\ \bar X^a(\si)\end{pmatrix}-\BR^{-1}(\Th^*\Bq+\Th^*_T\Bg)\\
\ns\ds\qq-\BR^{-1}(\BBL-\BR)\[-\BBL_\si^{-1}(I-\Pi_\si)(\Th^*\BQ,\Th^*_T\BG)\begin{pmatrix}\bar X_\si\\ X^a(\si)\end{pmatrix}-\BBL_\si^{-1}(I-\Pi_\si)(\Th^*\Bq+\Th^*_T\Bg)\]\\
\ns\ds\q=-\BR^{-1}\[I-(\BBL-\BR)\BBL_\si^{-1}(I-\Pi_\si)\]\[(\Th^*\BQ,\Th^*_T\BG)\begin{pmatrix}\bar X_\si\\ X^a(\si)\end{pmatrix}+(\Th^*\Bq+\Th^*_T
\Bg)\].\ea\ee
Setting $\si=t$, we obtain our conclusion, making use of \rf{LL}. \endpf

\ms

The above gives a {\it causal state feedback representation} for the open-loop optimal control in an abstract from. The appearance of $\BBL^{-1}$ makes the result hard to use since $\BBL$ is a complicated nonlocal operator. Hence, our next goal is to make the representation more explicitly accessible. We will achieve this goal in next section.

\section{Representation via Fredholm Integral Equations.}

Based on the result given in Theorem \ref{representation1}, we will now focus on the further manipulation of the abstract operator $\BBL^{-1}$. We want to convert it into another feedback gain operator which can be accessed in a computational manner. To this aim, for each $\si\in[0,T)$, we define $M_\si:[0,T]\times[0,T]\to\dbR^{m\times m}$ by the following:
\bel{M}\ba{ll}
\ds M_\si(t,s)=-R^{-1}(t)\Big(\[I-(\BBL-\BR)(I-\Pi_\si)\BBL^{-1}(I-\Pi_\si)\](\Th^*\BQ,\Th^*_T\BG)\begin{pmatrix} \Psi(\cd\,,s)\\ \Psi(T,s)\end{pmatrix}\Big)(t),\\
\ns\ds\qq\qq\qq\qq\qq\qq\qq\qq\qq\qq(t,s)\in [0,T]\times[0,T].\ea\ee
Note that $\Psi(\cd\,,\cd)$ is extended to be zero in $([0,T]\times[0,T])\setminus\Delta$. We would like to find an equation for $M_\sigma(\cd\,,\cd)$, which is easier to be used.

\bl{Lemma M} \sl For any $\si\in[0,T)$, the following Fredholm integral
equation
\bel{Fredholm}\ba{ll}
\ns\ds R(t)M_\si(t,s)=-\int_\si^T\(\int_{t\vee \xi}^T\Psi(\t,t)^\top Q(\t)
\Psi(\t,\xi)d\t+\Psi(T,t)^\top G\Psi(T,\xi)\)M_\si(\xi,s)d\xi
\\
\ns\ds\qq\qq\qq\q-\int_{t}^T\Psi(\t,t)^\top Q(\t)
\Psi(\t,s)d\t-\Psi(T,t)^\top G\Psi(T,s),\q (t,s)\in[0,T]\times[0,T],\ea\ee
admits a unique solution $M_\si(\cd,\cd)$, which is given by the expression in \rf{M}.

\el

\it Proof. \rm Note that (by \rf{LL}) and the fact that $\BR\Pi_\si=\Pi_\si\BR$,
\bel{L^{-1}}\ba{ll}
\ds(I-\Pi_\si)\BR^{-1}\big[I-(\BBL-\BR)(I-\Pi_\si)\BBL^{-1}(I-\Pi_\si)\big]\\
\ns\ds=(I-\Pi_\si)\BR^{-1}\big[(I-\Pi_\si)-(I-\Pi_\si)(\BBL-\BR)(I-\Pi_\si)\BBL^{-1}
(I-\Pi_\si)\big]\\
%
%\ns\ds=(I-\Pi_\si)\BR^{-1}\big[(I-\Pi_\si)-(I-\Pi_\si)\(\BBL(I-\Pi_t)\BBL^{-1}-\BR)\BBL^{-1}
%(I-\Pi_\si)\big]\\
%
%\ns\ds=(I-\Pi_\si)\BR^{-1}\big[(I-\Pi_\si)\BR\BBL^{-1}(I-\Pi_\si)\big]\\
%\ns\ds=(I-\Pi_\si)\BR^{-1}\big[I-(\BBL-\BR)\BBL^{-1}\big](I-\Pi_\si)\\
%
\ns\ds=(I-\Pi_\si)\BR^{-1}\big[(I-\Pi_\si)\BR(I-\Pi_\si)\BBL^{-1}(I-\Pi_\si)\big]
=(I-\Pi_\si)\BBL^{-1}(I-\Pi_\si).\ea\ee
Thus,
$$\ba{ll}
\ns\ds\q M_\si(t,s)=-R^{-1}(t)\Big\{\[I-(\BBL-\BR)(I-\Pi_\si)\BBL^{-1}(I-\Pi_\si)\](\Th^*\BQ,\Th^*_T\BG)\begin{pmatrix} \Psi(\cd\,,s)\\ \Psi(T,s)\end{pmatrix}\Big\}(t)\\
%
%\ns\ds=-R^{-1}(t)\Big(\Big[I-(\BBL-\BR)\(\BR^{-1}-\BR^{-1}(I-\Pi_\si)(\BBL-\BR)
%(I-\Pi_\si)\BBL^{-1}\)(I-\Pi_\si)\Big](\Th^*\BQ,\Th^*_T\BG)\begin{pmatrix} \Psi(\cd\,,s)\\ \Psi(T,s)\end{pmatrix}\Big)(t)\\
%
\ns\ds=-R^{-1}(t)\Big(\Big[I-(\BBL-\BR)(I-\Pi_\si)\BR^{-1}\(I-(\BBL-\BR)
(I-\Pi_\si)\BBL^{-1}
(I-\Pi_\si)\)\Big](\Th^*\BQ,\Th^*_T\BG)\begin{pmatrix} \Psi(\cd\,,s)\\ \Psi(T,s)\end{pmatrix}\Big)(t)\\
\ns\ds=-R^{-1}(t)\Big((\Th^*\BQ,\Th^*_T\BG)\begin{pmatrix} \Psi(\cd\,,s)\\ \Psi(T,s)\end{pmatrix}+(\BBL-\BR)(I-\Pi_\si)M_\si(\cd,s)\Big)(t)\\
\ns\ds=-R^{-1}(t)\Big(\big[\Th^*Q\Psi(\cd\,,s)\big](t)+\big[\Th^*_TG\Psi(T,s)\big](t)\Big)\\
\ns\ds\q-R^{-1}(t)\Big((\Th^*\BQ\Th+\Th_T^*\BG\Th_T)
(I-\Pi_\si)M_\si(\cd\,,s)\Big),\qq(t,s)\in[0,T]\times[0,T].\ea$$
This is equivalent to the following:
$$\ba{ll}
\ns\ds M_\si(t,s)=-R^{-1}(t)\Big(\int_t^T\Psi(\t,t)^\top Q(\t)\Psi(\t,s)d\t+\Psi(T,t)^\top G\Psi(T,s)\Big)\\
\ns\ds\qq\qq\q-R^{-1}(t)\int_t^T\Psi(\t,t)^\top Q(\t)\int_0^\t\Psi(\t,\xi)
{\bf1}_{[\si,T]}(\xi)M_\si(\xi,s)d\xi d\t\\
\ns\ds\qq\qq\q-R^{-1}(t)\Psi(T,t)^\top G\int_0^T\Psi(T,\xi){\bf1}_{[\si,T]}(\xi)M_\si(\xi,s)d\xi\\
\ns\ds\qq\qq=-R^{-1}(t)\int_\si^T\(\int_{t\vee \xi}^T\Psi(\t,t)^\top Q(\t)
\Psi(\t,\xi)d\t+\Psi(T,t)^\top G\Psi(T,\xi)\)M_\si(\xi,s)d\xi\\
\ns\ds\qq\qq\q-R^{-1}(t)\Big(\int_t^T\Psi(\t,t)^\top Q(\t)\Psi(\t,s)d\t+\Psi(T,t)^\top G\Psi(T,s)\Big),\qq(t,s)\in[0,T]\times[0,T].\ea$$
This means that $M_\si(\cd,\cd)$ is a solution to the Fredholm equation
\rf{Fredholm}. Now, for the uniqueness, it suffices to show that if
$$R(t)M_\si(t,s)+\int_\si^T\(\int_{t\vee \xi}^T\Psi(\t,t)^\top Q(\t)
\Psi(\t,\xi)d\t+\Psi(T,t)^\top G\Psi(T,\xi)\)M_\si(\xi,s)d\xi=0,$$
then $M_\si(t,s)=0$, $(t,s)\in[0,T]\times[0,T]$. Note that
$$\ba{ll}
\ns\ds 0=R(t)M_\si(t,s)+\Big((\BBL-\BR)(I-\Pi_\si)M_\si(\cd,s)\Big)(t)\\
\ns\ds\q=\Big(\Big[\BR\Pi_\si+\BBL(I-\Pi_\si)\Big]M_\si(\cd,s)\Big)(t),\q(t,s)\in[0,T]\times[0,T].\ea$$
Applying $(I-\Pi_\si)$ to the above, one has
$$0=\Big[(I-\Pi_\si)\BBL(I-\Pi_\si)M_\si(\cd,s)\Big](t)=\Big[\BBL_\si(I-\Pi_\si)M_\si(\cd,s)\Big](t),\q(t,s)\in[0,T]\times[0,T].$$
By the invertibility of $\BBL_\si$, one has $\Big((I-\Pi_\si)M_\si(\cd,s)\Big)(t)=0$, $(t,s)\in[0,T]\times[0,T]$. Consequently, $R(t)M_\si(t,s)=0,$ $(t,s)\in[0,T]\times[0,T]$. Hence, it follows from the invertibility of $\BR$ that
$M_\si(t,s)=0$, $(t,s)\in[0,T]\times[0,T]$ completing the proof. \endpf

\ms

We now prove the following theorem.

\bt{feedback} \sl Let {\rm(A5)} hold. Let $(\bar X(\cd),\bar u(\cd))$ be the open-loop optimal pair and $(\bar X_\si(\cd),\bar X^a(\cd))$ be the corresponding truncation and auxiliary trajectories. Then the open-loop optimal
control $\bar u(\cd)$ admits the following representation:
\bel{rep*}\ba{ll}
\ns\ds\bar u(t)=-R(t)^{-1}\(\big[\Th^*Q\bar X_t(\cd)\big](t)+\big[\Th^*_TG\bar X^a(t)\big](t)\)-R(t)^{-1}\(\big[\Th^*q(\cd)\big](t)+\big[\Th^*_Tg\big](t)\)
\\
\ns\ds\qq\q-\int_t^TM_t(t,s)R(s)^{-1}\(\big[\Th^*Q\bar X_t(\cd)\big](s)+\big[\Th^*_TG\bar X^a(t)\big](s)\)ds
\\
\ns\ds\qq\q
-\int_t^TM_t(t,s)R(s)^{-1}\(\big[\Th^*q(\cd)\big](s)+\big[\Th^*_Tg\big](s)\)ds,\qq\ae~t\in[0,T],\ea\ee
where $M_\si(\cd,\cd)$ is the unique solution of Fredholm equation \rf{Fredholm}.

\et

\it Proof. \rm First, similar to \rf{L^{-1}}, we have
\bel{L^{-1}*}\ba{ll}
\ds\[(I-\Pi_t)\BR^{-1}-(I-\Pi_t)\BBL^{-1}(I-\Pi_t)(\BBL-\BR)\BR^{-1}\](I-\Pi_t)\\
\ns\ds=\[(I-\Pi_t)-(I-\Pi_t)\BBL^{-1}(I-\Pi_t)\(\BBL(I-\Pi_t)-\BR\)\]\BR^{-1}(I-\Pi_t)\\
\ns\ds=\[(I-\Pi_t)-(I-\Pi_t)+(I-\Pi_t)\BBL^{-1}(I-\Pi_t)\BR\]
\BR^{-1}(I-\Pi_t)=(I-\Pi_t)\BBL^{-1}(I-\Pi_t).\ea\ee
Let us denote
$$\BBG(t)=(\Th^*\BQ,\Th^*_T\BG)\begin{pmatrix}\bar X_t\\ \bar X^a(t)\end{pmatrix}(t)+(\Th^*,\Th^*_T)\begin{pmatrix}q(\cd)\\ g\end{pmatrix}(t),\qq\ae~t\in[0,T].$$
Then,
$$\ba{ll}
\ns\ds\q\bar u(t)=-R(t)^{-1}\[I-(\BBL-\BR)(I-\Pi_t)\BBL^{-1}(I-\Pi_t)\]\BBG(t)\\
\ns\ds=-\BR^{-1}\BBG(t)+\BR^{-1}(\BBL-\BR)\Big[(I-\Pi_t)\BR^{-1}-(I-\Pi_t)\BBL^{-1}(I-\Pi_t)(\BBL-\BR)
\BR^{-1}\](I-\Pi_t)\BBG(t)\\
\ns\ds=-\BR^{-1}\BBG(t)+\BR^{-1}(\Th^*\BQ,\Th^*_T\BG)\begin{pmatrix}\Th\\ \Th_T\end{pmatrix}(I-\Pi_t)\BR^{-1}\BBG(t)\\
\ns\ds\q-\BR^{-1}(\BBL-\BR)(I-\Pi_t)\BBL^{-1}(I-\Pi_t)
(\Th^*\BQ,\Th^*_T\BG)\begin{pmatrix}\Th\\ \Th_T\end{pmatrix}(I-\Pi_t)\BR^{-1}\BBG(t)\\
\ns\ds=-\BR^{-1}\BBG(t)+\BR^{-1}\[I-(\BBL-\BR)(I-\Pi_t)\BBL^{-1}(I-\Pi_t)\]
(\Th^*\BQ,\Th^*_T\BG)\begin{pmatrix}\Th\\ \Th_T\end{pmatrix}(I-\Pi_t)\BR^{-1}\BBG(t)\ea$$
$$\ba{ll}
\ns\ds=-\BR^{-1}\[(\Th^*\BQ,\Th^*_T\BG)\begin{pmatrix}\bar X_t\\ \bar X^a(t)\end{pmatrix}(t)+(\Th^*,\Th^*_T)\begin{pmatrix}q(\cd)\\ g\end{pmatrix}(t)\]\\
\ns\ds\q+\BR^{-1}\[I-(\BBL-\BR)(I-\Pi_t)\BBL^{-1}(I-\Pi_t)\]
(\Th^*\BQ,\Th^*_T\BG)\\
\ns\ds\qq\cd\int^{T}_{t}\begin{pmatrix}\Psi(\cd,s)\\\Psi(T,s)\end{pmatrix} R(s)^{-1}\[\([\Th^*Q\bar X_t(\cd)
+\Th^*_TG\bar X^a(t)\)(s)+\big[\Th^*q(\cd)\big](s)
+\big[\Th^*_Tg\big](s)\]ds\\
\ns\ds=-R(t)^{-1}\(\big[\Th^*Q\bar X_t(\cd)\big](t)+\big[\Th^*_TG\bar X^a(t)\big](t)\)
-\int_t^TM_t(t,s)R(s)^{-1}\(\big[\Th^*Q\bar X_t(\cd)\big](s)+\big[\Th^*_TG\bar X^a(t)\big](s)\)ds
\\
\ns\ds\q-R(t)^{-1}\(\big[\Th^*q(\cd)\big](t)+\big[\Th^*_Tg\big](t)\)
-\int_t^TM_t(t,s)R(s)^{-1}\(\big[\Th^*q(\cd)\big](s)+\big[\Th^*_Tg\big](s)\)ds,\q\ae~t\in[0,T].\ea$$
This completes the proof. \endpf

\ms

Now, we return to the general case, i.e., $S(\cd)$ and $\rho(\cd)$ are not necessarily
zero. In this case, we summarize the result as follows:
\bel{hA}\left\{\2n\ba{ll}
\ds\h A(t,s)=A(t,s)-B(t,s)R(s)^{-1}S(s),\qq
\h\f(t)=\f(t)-\int_0^t{B(t,s)R(s)^{-1}\rho(s)\over(t-s)^{1-\b}}ds,\\
\ns\ds\h Q(s)=Q(s)-S(s)^\top R(s)^{-1}S(s),\qq
\h q(s)=q(s)-S(s)^\top R(s)^{-1}\rho(s),\\
\ns\ds\h\F(t,s)={\h A(t,s)\over(t-s)^{1-\b}}+\int_s^t{\h A(t,\t)\h\F(\t,s)
\over(t-\t)^{1-\b}}d\t,\\
\ns\ds\h\Psi(t,s)={B(t,s)\over(t-s)^{1-\b}}+\int_s^t{\h\F(t,\t)B(\t,s)\over
(\t-s)^{1-\b}}d\t,\qq\h\psi(t)=\h\f(t)+\int_0^t\h\F(t,s)\h\f(s)ds,\\
\ns\ds(\h\Th v)(t)=\int_0^t\h\Psi(t,s)v(s)ds,~ v\in\sU,~t\in[0,T],\q\h\Th_Tv=\int_0^T\h\Psi(T,s)v(s)ds=(\h\Th v)(T),\q v\in\sU.\ea\right.\ee
The truncation and auxiliary trajectories are defined by
\bel{auxiliary}\ba{ll}
\ns\ds X_\si(t)=\h\psi(t)+\int_0^{t\land\si}\h\Psi(t,s)\(u(s)+R(s)^{-1}\big[S(s)X(s)+\rho(s)\big]\)ds,\\
\ns\ds X^a(t)=\h\psi(T)+\int_0^t\h\Psi(T,s)\(u(s)+R(s)^{-1}\big[S(s)X(s)+\rho(s)\big]\)ds,\ea\qq t\in[0,T].\ee
The corresponding Fredholm equation reads
\bel{Fredholm*}\ba{ll}
\ns\ds R(t)\h M_\si(t,s)=-\int_\si^T\(\int_{t\vee \xi}^T\h\Psi(\t,t)^\top\h Q(\t)
\h\Psi(\t,\xi)d\t+\h\Psi(T,t)^\top G\h\Psi(T,\xi)\)\h M_\si(\xi,s)d\xi
\\
\ns\ds\qq\qq\qq\q-\int_{t}^T\h\Psi(\t,t)^\top\h Q(\t)
\h\Psi(\t,s)d\t+\h\Psi(T,t)^\top G\h\Psi(T,s),\q (t,s)\in[0,T]\times[0,T].\ea\ee
Then we can state the following result whose proof is clear.

\bt{rep-general} \sl Let {\rm(A2)--(A4)} hold. Let $(\bar X(\cd),\bar u(\cd))$ be the open-loop optimal pair of Problem (P). Then the open-loop
optimal control $\bar u(\cd)$ admits the following causal state feedback representation:
\bel{representation-general}\ba{ll}
\ns\ds\bar u(t)=-R(t)^{-1}\big[S(t)\bar X(t)+\rho(t)\big]-R(t)^{-1}\(\big[\h\Th^*\h Q\bar X_t(\cd)\big](t)+\big[\h\Th^*_TG\bar X^a(t)\big](t)\)\\
\ns\ds\qq\q-\int_t^T\h M_t(t,s)R(s)^{-1}\(\big[\h\Th^*\h Q\bar X_t(\cd)\big](s)+\big[\h\Th^*_TG\bar X^a(t)\big](s)\)ds
\\
\ns\ds\qq\q
-\int_t^T\h M_t(t,s)R(s)^{-1}\(\big[\h\Th^*\h q(\cd)\big](s)+\big[\h\Th^*_Tg\big](s)\)ds,
\\
\ns\ds\qq\q-R(t)^{-1}\(\big[\h\Th^*\h q(\cd)\big](t)+\big[\h\Th^*_Tg\big](t)\),\qq\ae~t\in[0,T],\ea\ee
where $\h\Th$, $\h\Th_T,$ $\h M_\si(\cd,\cd)$, $\h Q(\cd)$, and $\h q(\cd)$ are given
by \rf{hA} and \rf{Fredholm*}, $\bar{X}_\si(\cd)$ and $\bar{X}^a(\cd)$ are defined by \rf{auxiliary} with $u(\cd)$ being replaced by $\bar{u}(\cd)$.

\et

Since the general Volterra integral equation does not have a semigroup evolutionary property, the direct feedback implementation of the optimal control in terms of the actual trajectory $X(\tau),$ $\ae$ $\tau\in[0,T],$ is not possible, because the future information of the function $\psi$ is not counted. In view of this, the causal
stated feedback representation is the best that can be hoped for. Because
the auxiliary trajectory depends on $\h\psi(T)$, one might also call the above semi-causal state feedback representation as in \cite{Prichard-You 1996}.

\section{An Iteration Scheme for the Fredholm Integral Equation.}

In this section, we will briefly present a possible numerical scheme which is applicable to solve Fredholm integral equation (\ref{Fredholm}). This will, in principle, make the approach presented in the previous sections practically feasible.

\ms

During the period 1960--1990, there has been much work on developing and analyzing numerical methods for solving linear Fredholm integral equations of the second kind. The Galerkin and collocation methods are the well-established numerical methods (see \cite{Atkinson 1997}). Also, it is known that both the {\bf iterated} Galerkin and the {\bf iterated} collocation methods exhibit a higher order of convergence than the Galerkin method and collocation methods, respectively (see, for example \cite{Sloan 1984}). Long--Nelakanti \cite{Long-Nelakanti 2007} proposed an efficient iteration algorithm having much higher order of convergence, while they need less additional computational efforts for the implementation. Making use of the similar idea of \cite{Long-Nelakanti 2007}, we aim to obtain an efficient iteration scheme for the Fredholm integral equation (\ref{Fredholm}). Let us make this more precise now. To this end, we denote
$$\cK(t,\xi)=-R(t)^{-1}\(\int_{t\vee \xi}^T\Psi(\t,t)^\top Q(\t)
\Psi(\t,\xi)d\t+\Psi(T,t)^\top G\Psi(T,\xi)\),\q(t,\xi)\in[0,T]\times[0,T],$$
and
$$f(t,s)=-R(t)^{-1}\Big(\int_{t}^T\Psi(\t,t)^\top Q(\t)
\Psi(\t,s)d\t+\Psi(T,t)^\top G\Psi(T,s)\Big),\qq(t,s)\in[0,T]\times[0,T].$$
Next, we denote $\sM=L^2(0,T;\dbR^{m\times m})$, and for any $\sigma\in[0,T)$, define the integral operator $\sK_\si:\sM\rightarrow\sM$ by the following:
$$\sK_\si\eta(t)=\int_{\sigma}^T\cK(t,\xi)\eta(\xi)d\xi,\qq t\in[0,T],\q\forall \eta\in\sM.$$
Then, for each $\si\in[0,T)$, \rf{Fredholm} can be reconsidered as the following: for each $s\in[0,T]$,
$$(I-\sK_\si)M_\si(t,s)=f(t,s),\qq t\in[0,T].$$
We introduce a partition $\pi: 0=s_0<s_1<s_2<\cdots<s_N=T$ of $[0,T]$.
For each $s_i$, $0\les i\les N$, we consider the following equation:
\bel{Fredholm i}\ba{ll}
\ns\ds (I-\sK_\si)M_\si(t,s_i)=f(t,s_i),\qq t\in[0,T].\ea\ee
By assumption (A3), we can easily obtain that the integral operator $\sK_\si$ is a compact linear operator on $\sM$. Indeed,
$$\ba{ll}
\ns\ds\qq\int_0^T\int_0^T\Big|\cK(t,\xi)\Big|^2d\xi dt\\
\ns\ds\les K\int_0^T\int_0^T\Big|\int_{t\vee \xi}^T\Psi(\t,t)^\top Q(\t)
\Psi(\t,\xi)d\t\Big|^2d\xi dt+K\int_0^T\int_0^T\Big|\Psi(T,t)^\top G\Psi(T,\xi)\Big|^2d\xi dt\\
\ns\ds\les K\int_0^T\int_0^T\Big[\int_{t}^T|\Psi(\t,t)^\top|^2d\t\cdot\int_{\xi}^T|\Psi(\t,\xi)|^2d\t\Big] d\xi dt+K\int_0^T\int_0^T\frac{1}{(T-t)^{2(1-\b)}}\frac{1}{(T-\xi)^{2(1-\b)}}d\xi dt\\
\ns\ds\les  K\int_0^T\int_0^T\Big[\int_{t}^T\frac{1}{(\t-t)^{2(1-\b)}}d\t\cdot\int_{\xi}^T\frac{1}{(\t-\xi)^{2(1-\b)}}d\t\Big] d\xi dt
+K\int_0^T\frac{1}{(T-t)^{2(1-\b)}}dt<\infty,\ea$$
which implies that $\cK(\cd\,,\cd)\in L^2 ((0,T)\times(0,T);\dbR^{m\times m})$. Then, it is easy to see that $\sK_\si$ is compact on $\sM$ (see, for example Theorem 6.12 in \cite{Brezis}).
For each $s_i,$ $0\les i \les N,$ the Fredholm alternative theorem then guarantees the existence of a unique solution of (\ref{Fredholm i}) in $\sM$ (see, for example Theorem 1.3.1 in \cite{Atkinson 1997}). Let $\{\sM_n:n\ges 1\}$ be a sequence of increasing finite dimensional subspaces of $\sM$. Let $P_n: \sM\to\sM_n$ be the orthogonal projection operator (see, for example Section 3.1.2 in \cite{Atkinson 1997}).
Then, for each $s_i,$ $0\les i \les N,$ the Galerkin approximation is the solution of
\bel{Galerkin}M_{\si n}(t,s_i)=P_nf(t,s_i)+P_n\sK_\si M_{\si n}(t,s_i),\qq t\in[0,T].\ee
We can show that $P_ny\to y$, as $n\to\infty,$ for all $y\in\sM$ (see, for example Section 3.3.1 in \cite{Atkinson 1997}). Then, it follows from the compactness of $\sK_\si$ that
\bel{K convergence}\lim_{n\to\infty}\|\sK_\si-P_n\sK_\si\|_{\sL(\sM)}=0,\ee
(see, for example Lemma 3.1.2 in \cite{Atkinson 1997} or \cite{Sloan 1984}).
Then, by (\ref{K convergence}), one has that $(I-P_n\sK_\si)^{-1}$ exist and uniformly bounded for sufficiently large $n,$ and the approximation scheme is uniquely solvable (see, for example Theorem 3.1.1 in \cite{Atkinson 1997}).

For each $s_i,$ $0\les i \les N,$ the iterated Galerkin approximation (see, \cite{Sloan 1984}) may be defined by
\bel{iterated Galerkin}\widetilde{M}_{\si n}(t,s_i)=f(t,s_i)+\sK_\si M_{\si n}(t,s_i),\qq t\in[0,T].\ee
Applying $P_n$ to both side of (\ref{iterated Galerkin}), we have $P_n\widetilde{M}_{\si n}=M_{\si n}$, and hence for each $s_i,$ $0\les i \les N,$
it holds that
$$\widetilde{M}_{\si n}(t,s_i)=f(t,s_i)+\sK_\si P_n \wt M_{\si n}(t,s_i),\qq t\in[0,T].$$
One can show that the iterated Galerkin scheme (\ref{iterated Galerkin}) can converge more rapidly than the rate achieved by the approximation (\ref{Galerkin}) (see, \cite{{Graham-Joe-Sloan 1985}} and \cite{Sloan 1984}). Further, Long--Nelakanti \cite{Long-Nelakanti 2007} proposed a more efficient iteration algorithm which even has much higher order of convergence than the iterated Galerkin scheme. The iteration algorithm is as follows: for each $s_i,$ $0\les i \les N,$
set $M^{(0)}_{\si n}(t,s_i)=\widetilde{M}_{\si n}(t,s_i),$ $ t\in[0,T],$ then for $k=0,1,\ldots,$
$$\ba{ll}
\ns\ds\mbox{step}~1:\q \wt M^{(k)}_{\si n}(t,s_i)=f(t,s_i)+\sK_\si M^{(k)}_{\si n}(t,s_i),\qq t\in[0,T],\\
\ns\ds\mbox{step}~2:\q \wt{\wt{M}}^{(k)}_{\si n}(t,s_i)=f(t,s_i)+\sK_\si \wt{M}^{(k)}_{\si n}(t,s_i),\qq t\in[0,T],\\
\ns\ds\mbox{step}~3:\q g^{(k)}_{n}(t,s_i)=\widetilde{\widetilde{M}}^{(k)}_{\si n}(t,s_i)-\widetilde{M}^{(k)}_{\si n}(t,s_i),\qq t\in[0,T],\\
\ns\ds\mbox{step}~4:\q \mbox{for}~\mbox{each}~s_i,~ 0\les i \les N,~ \mbox{seeking}~\mbox{a}~\mbox{unknown}~\mbox{function}~e^{(k)}_{n}(t,s_i),~ t\in[0,T]~\mbox{by}~\mbox{solving}~\mbox{the}~\mbox{equation}\\
\ns\ds\qq\qq~
(I-P_n\sK_\si)e^{(k)}_{n}(t,s_i)=P_ng^{(k)}_{n}(t,s_i),\qq t\in[0,T],\\
\ns\ds\mbox{step}~5:\q M^{(k+1)}_{\si n}(t,s_i)=\sK_\si e^{(k)}_{n}(t,s_i)+\widetilde{\widetilde{M}}^{(k)}_{\si n}(t,s_i),\qq t\in[0,T].\ea$$
By the superconvergence rates for every step of iteration, we obtain that for each  $k=0,1,\ldots,$ for each $s_i,$ $0\les i <N,$
$$\ba{ll}
\ns\ds M^{(k+1)}_{\si n}(\cd,s_{i})\to M_{\si}(\cd,s_{i})~\mbox{in}~\sM,~\mbox{as}~n\to\infty.\ea$$
For more details, see \cite{Long-Nelakanti 2007}.
Thus, on a subsequence, still denoted in the same way, for each $n=0,1,\ldots,$
\bel{superconvergence}\|M^{(k+1)}_{\si n}(\cd,s_{i})-M_{\si}(\cd,s_{i})\|_{\sM}\les\frac{1}{n}.\ee
Now, for $k=0,1,\ldots,$ let
$$M^{N}_{\si n }(t,s)=\sum^{N}_{i=1}\Big[\frac{s_i-s}{s_i-s_{i-1}}M^{(k+1)}_{\si n}(t,s_{i-1})+\frac{s-s_{i-1}}{s_i-s_{i-1}}M^{(k+1)}_{\si n}(t,s_i)\Big]{\bf1}_{[s_{i-1},s_i)}(s),\q (t,s)\in[0,T]\times[0,T].$$
Then, for $k=0,1,\ldots,$
$$\ba{ll}
\ns\ds\q|M^{N}_{\si n }(t,s)-M_{\si }(t,s)|\\
\ns\ds\les\sum^{N}_{i=1}\Big[\frac{s_i-s}{s_i-s_{i-1}}\Big(|M^{(k+1)}_{\si n}(t,s_{i-1})-M_{\si}(t,s_{i-1})|+|M_{\si}(t,s_{i-1})-M_{\si}(t,s)|\Big)\\
\ns\ds\qq\q+\frac{s-s_{i-1}}{s_i-s_{i-1}}\Big(|M^{(k+1)}_{\si n}(t,s_{i})-M_{\si}(t,s_{i})|+|M_{\si}(t,s_{i})-M_{\si}(t,s)|\Big)\Big]{\bf1}_{[s_{i-1},s_i)}(s)\ea$$
$$\ba{ll}
\ns\ds\les \sum^{N}_{i=1}\Big[\frac{s_i-s}{s_i-s_{i-1}}|M_{\si}(t,s_{i-1})-M_{\si}(t,s)|+\frac{s-s_{i-1}}{s_i-s_{i-1}}|M_{\si}(t,s_{i})-M_{\si}(t,s)|\Big]{\bf1}_{[s_{i-1},s_i)}(s)\\
\ns\ds\q+K\sum^{N}_{i=1}\Big[|M^{(k+1)}_{\si n}(t,s_{i-1})-M_{\si}(t,s_{i-1})|+|M^{(k+1)}_{\si n}(t,s_{i})-M_{\si}(t,s_{i})|\Big]{\bf1}_{[s_{i-1},s_i)}(s),\q (t,s)\in[0,T]\times[0,T].\ea$$
Hence, for $k=0,1,\ldots,$ for each $s\in[0,T),$
\bel{Mcontinuous}\ba{ll}
\ns\ds\q \|M^{N}_{\si n}(\cd,s)-M_{\si }(\cd,s)\|_{\sM}\\
\ns\ds\les K\sum^{N}_{i=1}\Big[\frac{s_i-s}{s_i-s_{i-1}}\|M_{\si}(\cd,s_{i-1})-M_{\si}(\cd,s)\|_{\sM}+\frac{s-s_{i-1}}{s_i-s_{i-1}}\|M_{\si}(\cd,s_{i})-M_{\si}(\cd,s)\|_{\sM}\Big]{\bf1}_{[s_{i-1},s_i)}(s)\\
\ns\ds\q+ K\sum^{N}_{i=1}\Big[\|M^{(k+1)}_{\si n}(\cd,s_{i-1})-M_{\si}(\cd,s_{i-1})\|_{\sM}+\|M^{(k+1)}_{\si n}(\cd,s_{i})-M_{\si}(\cd,s_{i})\|_{\sM}\Big]{\bf1}_{[s_{i-1},s_i)}(s).\ea\ee
We need a little more assumption which we now introduce.

\ms

{\bf(A6)} \rm There exists a modulus of continuity
$\o(\cd)$ such that
$$|B(t,s')-B(t,s)|\les\o(s'-s),\qq (t,s'),~(t,s)\in\D.$$
Then, pick any $s_0\in[0,T),$  making use of the similar argument of the proof for Corollary \ref{corollary 2.2}, (ii), we obtain that
$\Psi(\cd,s)$ is continuous in $L^2(0,T;\dbR^{n\times m})$ at $s_0,$ and $\Psi(T,s)$ is continuous at $s_0.$
Since $-\BR^{-1}\[I-(\BBL-\BR)\BBL_\si^{-1}(I-\Pi_\si)\]\Th^*\BQ$ is a bounded linear operator from $L^2(0,T;\dbR^{n\times m})$ to $\sM$ and $-\BR^{-1}\[I-(\BBL-\BR)\BBL_\si^{-1}(I-\Pi_\si)\]\Th^*_T\BG$ is a bounded linear operator from $\dbR^{n\times m}$ to $\sM$, it holds that $M_{\si}(\cd\,,s)$ is continuous in $\sM$ at $s_0\in[0,T).$

\ms

Consequently, for $k=0,1,\ldots$, for each $s\in[0,T)$, we let $N\rightarrow\infty$ and $n\rightarrow\infty.$ Then, by (\ref{superconvergence}) and (\ref{Mcontinuous}), we obtain that $M^{N}_{\si n}(\cd\,,s)\to M_{\si }(\cd\,,s)$ in $\sM$. Thus, we obtain a feasible numerical scheme which can be used to solve Fredholm integral equation (\ref{Fredholm}).

\section{Concluding Remarks.}

In this paper, we have studied an optimal control problem, with the state equation being a linear Volterra integral equation having a singular kernel. The cost functional is a form of quadratic plus linear terms of the state and the control. Under proper conditions, the open-loop optimal control uniquely exists. However, normally, the open-loop optimal control is not of non-anticipating form. Our main goal is to obtain a causal state feedback representation of the open-loop optimal control. In doing that, we have introduced a Fredholm integral equation which plays a role of Riccati equation in the standard LQ problems for ODE systems. To make our result practically feasible (in principle), we have briefly presented a possible numerical scheme for computing the solution to the Fredholm integral equation.

\end{document}